\documentclass[leqno,twoside]{article}
\usepackage{amsmath,amssymb,amsthm}
\usepackage[margin=2.5cm]{geometry} 

\usepackage{titlesec,hyperref}
\usepackage{fancyhdr}

\usepackage{color}

\fancypagestyle{plain}
{
\fancyhf{}

}

\titleformat{\subsection}{\it}{\thesubsection.\enspace}{1pt}{}

\newtheorem{theo}{Theorem}[section]
\newtheorem{lemm}[theo]{Lemma}
\newtheorem{defi}[theo]{Definition}

\newtheorem{prop}[theo]{Proposition}
\newtheorem{rema}[theo]{Remark}
\numberwithin{equation}{section}

\allowdisplaybreaks

\begin{document}
\title{Gevrey Smoothing Effect of Solutions to Non-Cutoff Boltzmann Equation for Soft Potential with Mild and Critical Singularity}

\author{Teng-Fei Zhang$^1$
\quad Zhaoyang Yin$^2$ \\[10pt]
Department of Mathematics, Sun Yat-sen University,\\
510275, Guangzhou, P. R. China.
}
\footnotetext[1]{Corresponding author. Email: {\it fgeyirui@163.com.} TEL.: {\it 18602078196.}}
\footnotetext[2]{Email: \it mcsyzy@mail.sysu.com.cn}
\date{}
\maketitle
\hrule

\begin{abstract}

In this paper we study the Gevrey smoothing effect of solutions to the non-cutoff spatially homogeneous Boltzmann equation for soft potential. We consider not only the mild singularity case $s<1/2$ as we did in the previous works for spatially homogeneous case (J. Differential Equations 253(4) (2012), 1172-1190. DOI: 10.1016/j.jde.2012.04.023) and for spatially inhomogeneous case ({Preprint \href{http://arxiv.org/abs/1304.2971v2}{arXiv: 1304.2971v2.pdf}}), but also the critical singularity case $s=1/2$ (with a particular soft potential $\gamma=-2$). Besides, we try to extend the range of $\gamma$. We derive a new coercivity estimate for collision operator, from which we can obtain the propagation of Gevrey regularity for $\gamma \in (-5/2,0)$, and Gevrey regularity for $\gamma \in [-2,0)$, which improve the previous assumption $\gamma \in (-1-2s,0)$. In addition, we consider $\gamma$ and $s$ separately instead of viewing $\gamma+2s$ as one untied quantity.

\vspace*{5pt}

\vspace*{5pt}
\end{abstract}

\vspace*{10pt}


\tableofcontents
\section{Introduction}

\subsection{The Boltzmann equation}

In this paper we consider the Cauchy problem of the non-cutoff Boltzmann equation. First we introduce
the Cauchy problem of the full (or, spatially inhomogeneous) Boltzmann equation without
angular cutoff, with a $T >0$,
\begin{align}\label{BE-full}
  \left\{
    \begin{array}{l}\displaystyle
     f_t(t,x,v)+v\cdot \nabla_xf(t,x,v)=Q(f,f),\quad t\in (0,T], ~x\in \mathbb{T}^3,~v \in \mathbb{R}^3,\\
     f(0,x,v)=f_0(x,v).
    \end{array}
  \right.
\end{align}

Above, $f=f(t,x,v)$ describes the density distribution function of particles located around position $x\in \mathbb{T}^3$ with velocity $v\in \mathbb{R}^3$ at time $t \geq 0$. The right-hand side of the first equation is the so-called Boltzmann bilinear collision operator acting only on the velocity
variable $v$:
\[
Q(g, f)=\int_{\mathbb{R}^3}\int_{\mathbb S^{2}}B\left({v-v_*},\sigma
\right)
 \left\{g'_* f'-g_*f\right\}d\sigma dv_*.
\]
Note that we use the well-known shorthand $f=f(t,x,v)$, $f_*=f(t,x,v_*) $, $f'=f(t,x,v') $, $f'_*=f(t,x,v'_*) $ throughout this paper.

In this paper we consider the Cauchy problem of the Boltzmann equation in the spatially homogeneous case, that is, for a $T >0$,
\begin{align}\label{BE}
  \left\{
    \begin{array}{l}
      f_t(t,v)=Q(f,f)(v),\quad t\in (0,T], ~  v \in \mathbb{R}^3,\\
      f(0,v)=f_0(v),
    \end{array}
  \right.
\end{align}
where ``spatially homogeneous" means that $f$ depends only on $t$ and $v$.

By using the $\sigma$-representation, we can describe the
relations between the post- and pre-collisional velocities as
follows, for $\sigma \in \mathbb S^2$,
$$
v'=\frac{v+v_*}{2}+\frac{|v-v_*|}{2}\sigma,~ ~ v'_*
=\frac{v+v_*}{2}-\frac{|v-v_*|}{2}\sigma.
$$
We point out that the collision process satisfies the conservation of momentum and kinetic energy, i.e.
$$
v+v_*=v'+v'_*,\qquad  |v|^2+|v_*|^2=|v'|^2+|v'_*|^2.
$$

The collision cross section $B(z, \sigma)$ is a given non-negative function depending only on the interaction law between particles. From a mathematical viewpoint, that means $B(z, \sigma)$ depends only on the relative velocity $|z|=|v-v_*|$ and the deviation angle $\theta$ through the scalar product
$\cos \theta=\frac{z}{|z|} \cdot \sigma$.

The cross section $B$  is assumed here to be of the type:
$$
B(v-v_*, \cos \theta)=\Phi (|v-v_*|) b(\cos \theta),~~
\cos \theta=\frac{v-v_*}{|v-v_*|} \cdot \sigma,~~
0\leq\theta\leq\frac{\pi}{2}.
$$
Above, $\Phi$ stands for the kinetic factor which is of the form:
$$
\Phi (|v-v_*|) = |v-v_*|^{\gamma},
$$
and $b$ denotes the angular part with a singularity such that,
$$
\sin \theta b(\cos \theta) \sim K\theta^{-1-2s}, \ \ \mbox{as} \ \ \theta\rightarrow 0+,
$$
for some positive constant $K$ and $0< s <1$.

We remark that if the inter-molecule potential satisfies the inverse-power law $U(\rho) = \rho ^{-(p-1)} ~(\textrm{where } p>2) $, it holds $\gamma = \frac{p-5}{p-1} $, $ s=\frac{1}{p-1} $. Generally, the cases $\gamma >0$, $\gamma =0$, and $\gamma <0$ correspond to so-called hard, Maxwellian, and soft potential respectively. And the cases $0<s<1/2$, $1/2 \leq s<1$ correspond to so-called mild singularity and strong singularity respectively.

As is discussed by Desvillettes in \cite{Desvillettes-2001} (also, by Villani in his handbook \cite{Villani}), the most interesting assumption for $\gamma$ is included in $(-3,1)$. Further, we announce that we will focus on the mild singularity case $0<s<1/2$ for soft potential $\gamma \in (-5/2,0)$.

\subsection{Review of related references}

Now we give a brief review about some related researches. Firstly we refer the reader to Villani's review book \cite{Villani} for the physical background and the mathematical theories of the Boltzmann equation. And for more details about the non-cutoff theories, we refer to Alexandre's review paper \cite{Alex-review}.

Before continuing the statement, we introduce the definition of Gevrey spaces $G^s(\Omega)$ where $\Omega$ is an open subset of $\mathbb{R}^3$. (It could be found in many references, e.g. \cite{Mori-Xu-09,Zhang-Yin}.)
\begin{defi}\label{Def_2}
For $0<s<+\infty$, we say that $ f \in G^s(\Omega) $, if $f \in C^\infty(\Omega)$, and there exist $C>0,~ N_0>0$ such that
$$ \|\partial^\alpha f\|_{L^2(\Omega)} \leq C^{|\alpha|+1} {\{\alpha! \}^s},\quad
\forall \alpha \in \mathbb{N}^3,~ |\alpha| \geq N_0.
$$

If the boundary of ~$\Omega$ is smooth, by using the Sobolev embedding theorem, we have the same type estimate with $L^2$ norm replaced by any $L^p$ norm for $2 < p \leq +\infty$.

When $s = 1$, it is usual analytic function. If $s > 1$, it is Gevrey class function. And for $0 < s < 1$, it is called ultra-analytic function.
\end{defi}

In 1984 Ukai showed in \cite{Ukai-84} that there exists a unique local solution to the Cauchy problem for the Boltzmann equation in Gevrey classes for both spatially homogeneous and inhomogeneous cases, under the assumption on the cross section:
\begin{align*}
&\big| B(|z|,\cos \theta) \big| \leq  K(1+|z|^{-\gamma'}+|z|^\gamma) \theta^{-n+1-2s},
      \quad n  \textrm{ is dimensionality},\\
&(0\leq \gamma' < n,~ 0\leq \gamma <2,~ 0\leq s<1/2, ~\gamma +6s<2 ).
\end{align*}

By introducing the norm of Gevrey space
$$
\|f\|^{U}_{\delta,\rho,\nu} = \sum_{\alpha} \frac{\rho^{|\alpha|}}{\{\alpha!\}^\nu}
\|e^{\delta \langle v \rangle^2} \partial_v^{\alpha} f \|_{L^\infty(\mathbb{R}^n_v)},
$$
Ukai proved that in the spatially homogeneous case, for instance, under some assumptions for $\nu$ and the initial datum $f_0(v)$, the Cauchy problem (\ref{BE}) has a unique solution $f(t,v)$ for $t\in (0,T]$.

In \cite{Desvillettes2} Desvillettes studied firstly the $C^\infty$ smoothing effect for solutions of Cauchy problem in spatially homogeneous non-cutoff case, and conjectured Gevrey smoothing effect. And he proved in \cite{Desvillettes} the propagation of Gevrey regularity for solutions without any assumptions on the decay at infinity in $v$ variables.

In 2009 Morimoto et al. considered in \cite{Mori-Ukai-Xu-Yang} the Gevrey regularity for the linearized Boltzmann equation around the absolute Maxwellian distribution, by virtue of the following mollifier:
$$
G_\delta (t,D_v)=\frac{e^{t \langle D_v \rangle^{1/\nu}}}
                      {1+\delta e^{t \langle D_v \rangle^{1/\nu}}},\quad 0< \delta <1.
$$
We remark that the same operator was used in many related researches and models such as the Kac's equation (a simplification of Boltzmann equation to one dimension case), the ultra-analytic smoothing effect for spatially homogeneous nonlinear Landau equation and the linear and non-linear Fokker-Planck equations.

In the mild singularity case $0<s<1/2$, Huo et al. proved in \cite{Huo} that any weak solution $f(t,v)$ to the Cauchy problem (\ref{BE}) satisfying the natural boundedness on mass, energy and entropy, namely,
\begin{align}\label{natural bound}
\int_{\mathbb{R}^n} f(v)[1+|v|^2+\log(1+f(v))]dv < +\infty,
\end{align}
belongs to $H^{+\infty}(\mathbb{R}^n)$ for any $0<t \leq T$, and moreover,
\begin{align}\label{smooth solution}
f \in L^\infty \big([t_0,T];H^{+\infty}(\mathbb{R}^n)\big),
\end{align}
for any $T>0$ and $t_0 \in (0,T)$.

In the framework of small perturbation of an equilibrium state, Alexandre et al. studied the Cauchy problem of the Boltzmann equation for soft and hard potential (see \cite{FiveGroup-I, FiveGroup-II}), and obtain the global existence of solution in weighted Sobolev spaces. Some other results about existence of perturbative solutions for the cutoff Boltzmann equation are due to Guo (for instance, see \cite{Guo}) around a global Maxwellian, and Tai-Ping Liu (compare \cite{Liu1, Liu2}, for example), and so on.

In \cite{FiveGroup-regulariz} the five authors considered a kind of solution having the Maxwellian decay, based on which we introduce the following definition:
\begin{defi}\label{Def1}
We say that $f(t,v)$ is a smooth Maxwellian decay solution to the Cauchy problem (\ref{BE}) if
\begin{align*}
  \left\{
     \begin{array}{l}
       f \geq 0,~ \not \equiv 0,\\
       \exists ~ \delta_0 >0 \textrm{ such that }
          e^{\delta_0 \langle v \rangle ^2} f\in L^\infty \left( [0,T];H^{+\infty}(\mathbb{R}^3) \right).
     \end{array}
  \right.
\end{align*}
\end{defi}

(Note that the Theorem 1.2 of \cite{FiveGroup-regulariz} shows the uniqueness of the smooth Maxwellian decay solution to the Cauchy problem (\ref{BE}).)

The five authors also proved in \cite{FiveGroup-regulariz} the smoothing effect on the solutions with weight. In detail, if the non-negative $f$ belongs to $\mathcal{H}_l^5 \Big( (t_1,t_2) \times \Omega \times \mathbb{R}^3_v \Big)$, solves the spatially inhomogeneous Boltzmann equation (\ref{BE-full}) in this domain in the classic sense, and satisfies the non-vacuum condition $\|f(t,x,v)\|_{L^1(\mathbb{R}^3_v)} >0$, then
$$
f\in \mathcal{H}_l^{+\infty} \Big( (t_1,t_2) \times \Omega \times \mathbb{R}^3_v \Big),
$$
hence it follows that,
$$
f\in C^\infty \Big( (t_1,t_2) \times \Omega; \mathcal{S}(\mathbb{R}^3_v) \Big).
$$

In 2010 Morimoto-Ukai considered in \cite{Mori-Ukai} the Gevrey regularity of $C^\infty$ solutions with the Maxwellian decay to the Cauchy problem of spatially homogeneous Boltzmann equation. Motivated by their idea, we considered the problem in \cite{Zhang-Yin} for a more general case. More precisely, we considered the general kinetic factor $ \Phi (|v|) = |v|^{\gamma} $ instead of the moderate form $ \langle v \rangle ^{\gamma}=( 1+|v|^2)^{\gamma/2} $ in \cite{Mori-Ukai}, and a wider range of the parameter of $\gamma$ (s.t. $\gamma +2s \in (-1,1)$) so as to fit for both hard potential and soft potential.

In the ensuing paper \cite{Zhang-Yin-2} we studied still the general case $\Phi (|v|) = |v|^{\gamma}$ with $\gamma +2s \in (-1,1)$ in the mild singularity case, but the spatially inhomogeneous case. We obtain a corresponding result about the Gevrey regularity.

In this present paper, we resume to the spatially homogeneous case, and try to extend the range of $\gamma$ for soft potential (in $\gamma \geq 0$ we make no change). For the propagation of Gevrey smoothing effect, we consider the case $\gamma \in (-5/2,0)$ here to take place of the previous assumption $\gamma \in (-1-2s,0)$. ({\it Recall that in \cite{Zhang-Yin}, we assume $\gamma+2s \in (-1,1)$. Note that the simple inequality $-1-2s>-2>-5/2$ gives the range of extending.}) To prove the order of Gevrey regularity, we assume further that $\gamma \in [-2,0)$ (so as to make the extra weight $\gamma/2$ to be no more than 1). We emphasize that we consider not only the mild singularity case $s<1/2$, but also the critical singularity case $s=1/2$ (with a particular soft potential $\gamma=-2$). It is generally known that there are few results concerning the strong singularity case $s \geq 1/2$, even if for the Maxwellian case. In addition, $\gamma$ and $s$ are considered separately instead of viewing $\gamma+2s$ as one united quantity. To achieve the goal, we need a new coercivity estimate for collision operator different from that the authors proved in \cite{Chen-He-1} with respect to the parameter $\gamma$.

\subsection{Main results}

Now we give our main result of propagation of Gevrey regularity in spatially homogeneous case as follows:
\begin{theo}\label{Propagation of Gevrey}
Let $\nu >1 $(which is independent of s) and assume that $0<s \le 1/2 $, $-5/2<\gamma<0$. Let $f(t,v)$ be a smooth Maxwellian decay solution to the Cauchy problem (\ref{BE}). If there exist $\rho'$, $\delta'$ such that
\begin{align}\label{initial assumption}
\sup_{\alpha} \frac{\rho'^{|\alpha|} \|e^{\delta' \langle v \rangle ^2}
              \partial^{\alpha}_{v} f(0)\|_{L^2}}{\{\alpha!\}^\nu} < +\infty,
\end{align}
then there exist $ \rho >0 $ and $ \delta,\kappa >0 $ with $ \delta >\kappa T $ such that
\begin{align}\label{propag}
\sup_{t\in (0,T]}\sup_{\alpha} \frac{\rho^{|\alpha|} \|e^{(\delta-\kappa t) \langle v \rangle^2}
                               \partial^{\alpha}_{v} f(t)\|_{L^2}}{\{\alpha!\}^\nu} < +\infty.
\end{align}
\end{theo}

\begin{rema}
It should be noted that the above theorem is similar as Theorem 1.2 in \cite{Mori-Ukai} and Theorem 1.3 in \cite{Zhang-Yin}, but we consider here $\Phi=|v|^\gamma$ and $\gamma \in (-5/2,0)$.
\end{rema}

Arguing as in Section 4 of \cite{Mori-Ukai}, and thanks to the interpolation inequality with weight (by Lemma 3.7 of \cite{FiveGroup-Smoothing}),
$$
  \|W_l F^{(\alpha,\beta)} \|^2_{H^{s/2}}
\lesssim  \|W_{l+\gamma/2} F^{(\alpha,\beta)} \|_{H^s}
          \|W_{l-\gamma/2} F^{(\alpha,\beta)} \|_{L^2}
\lesssim  \|W_{l+\gamma/2} F^{(\alpha,\beta)} \|_{H^s}
          \|W_{l+1} F^{(\alpha,\beta)} \|_{L^2},  \text{ for } \gamma \in [-2,0),
$$
we would obtain the Gevrey smoothing effect of order $1/s$ as follows:
\begin{theo}\label{Gevrey Regularity}
Assume that $ 0<s \le 1/2 $, $\gamma \in [-2,0)$. Let $ \nu =1/s$ and let $f(t,v)$ be a smooth Maxwellian decay solution to the Cauchy problem (\ref{BE}), then for any $ t_0 \in (0,T)$, there exist $ \rho >0 $ and $ \delta,\kappa >0 $ with $ \delta >\kappa T $ such that
\begin{align}
  \sup_{t\in [t_0,T]} \sup_{\alpha} \frac{\rho^{|\alpha|} \|e^{(\delta-\kappa t) \langle v \rangle^2}
                                    \partial^{\alpha}_{v} f(t)\|_{L^2}}{\{\alpha!\}^\nu} < +\infty.
\end{align}
\end{theo}

\begin{rema}
From the above interpolation inequality, we know that the additional requirement $\gamma \ge -2$ is natural, which is different from the assumption $\gamma \in (-5/2,0)$ in Theorem \ref{Propagation of Gevrey}.
\end{rema}

\begin{rema}
We believe that these results are also valid for spatially inhomogeneous case, as we did in \cite{Zhang-Yin-2}.
\end{rema}
\subsection{The structure of the paper}

The remainder of the paper proceeds as follows. In the next section we give some preliminaries and a main lemma, by using which we can complete the proof of Theorem \ref{Propagation of Gevrey} immediately. In Section 3 we prove a coercivity estimate for collision operator which is different from what we used before. The proof of the main lemma will be given in Section 4.

\section{Preliminaries}

First of all, we introduce some basic definitions (see \cite{Mori-Ukai} for details).

Let $l,~r \in \mathbb{Z}_+$ which will be chosen later. For $\delta,~ \rho >0$ we set:
$$
\|f\|_{\delta,l,\rho,\alpha,r}
\triangleq \frac{ \rho^{|\alpha|} \| \langle v \rangle ^l e^{\delta \langle v \rangle ^2} \partial^\alpha_v f\|_{L^2}}{\{(\alpha - r)!\}^\nu},
$$
where $ \alpha =(\alpha_1,\alpha_2,\ldots,\alpha_n) \in \mathbb{Z}^n_+ $, and we denote
$$ (\alpha - r)!=(\alpha_1 - r)! \cdots (\alpha_n - r)!.$$

    Now we give the following definition:
\begin{align}
\|f\|_{l,\rho,r,N}(t)
\triangleq \sup_{rn \leq |\alpha| \leq N} \|f\|_{\delta -\kappa t,l,\rho,\alpha,r},
\end{align}
with fixed $\delta,~ \kappa >0$ such that $\delta > \kappa T$. Here $N$  is a fixed large number satisfying $rn \leq |\alpha| \leq N$.
Then for $ h > 1$ we can obtain
\begin{align}
\|f\|_{l,\rho (1+h)^{-\nu},r,N}(t) \leq \left\{ \frac{(r!)^n}{h^r}\right\}^\nu \|f\|_{l,\rho,0,N}(t).
\end{align}

Now let $\rho =\rho'$ in the above inequality and take a large enough $ h $, then it follows from the initial assumption (\ref{initial assumption}) that $ \|f\|_{l,\rho' (1+h)^{-\nu},r,N}(0) $ is as small as we want, where $\delta$ can be chosen any positive less than $\delta' > 0$ in (\ref{initial assumption}).

Thus, to prove (\ref{propag}), it suffices to prove, under the assumption that $ \|f\|_{l,\rho,r,N}(0) $ is sufficiently small,
\begin{align}
\sup_{t \in (0,T]} \|f\|_{l,\rho,r,N}(t) < \infty.
\end{align}
Above, $\rho =\rho' (1+h)^{-\nu}$.

    We point out that we will consider the Cauchy problem in $\mathbb{R}^3$ in the paper.
\begin{lemm}\label{Main Lemma}
If $ l \geq 4$ and $r >1+\nu/(\nu-1)$ then for any $\alpha$ satisfying $ 3r \leq |\alpha| \leq N $ we have
\begin{align}\label{basic inequa}
&\|f(t)\|^2_{\delta -\kappa t,l,\rho,\alpha,r}
+ 2\kappa \int_0^t \|f(\tau)\|^2_{\delta-\kappa \tau,l+1,\rho,\alpha,r}d\tau \\
\leq &   \|f(0)\|^2_{\delta,l,\rho,\alpha,r}
       + C_\kappa \int_0^t \Big( \|f\|^2_{l,\rho,r,N}(\tau)
                             + \|f\|^{2(i+\beta)/\beta}_{l,\rho,r,N}(\tau) \Big)d\tau
       + \frac{\kappa}{10}\sup_{3r \leq |\alpha| \leq N}
                         \int_0^t \|f(\tau)\|^2_{\delta-\kappa \tau,l+1,\rho,\alpha,r}d\tau, \nonumber
\end{align}
where $i=\Large \emph {1}_{s<1/2} + 2 \cdot \Large \emph {1}_{s=1/2}$, and $\beta = -\gamma \cdot \Large \emph {1}_{\gamma \in [-1,0)} + \Large \emph {1}_{\gamma \in (-5/2,-1)}$.
\end{lemm}

Mimicking the scheme in Section 2 in \cite{Mori-Ukai} (It is only a matter of using the Bernoulli equation in ODE theory with a different order), we can prove Theorem \ref{Propagation of Gevrey} by virtue of this lemma. We omit the details. The proof of this lemma will be given in Section 4.

\section{Coercivity estimates}

Considering the coercivity estimate of collision operator $-Q(g,f)$ for $\gamma \in (-5/2,0)$, the result in \cite{Chen-He-1} for the case $\gamma+2s >-1$ is no longer applied to our meet. We reconsider in the framework of \cite{Chen-He-1} the coercivity estimate, which is described as follows:
\begin{lemm}\label{coercivity lemma}
Let $0<s<1$ and assume that the nonnegative function $g$ satisfies
$$
\|g\|_{L^1_2(\mathbb{R}^3_v)}+\|g\|_{L\log L(\mathbb{R}^3_v)}<\infty.
$$
Then there exists a constant $C_g>0$ depending on $B$, $\|g\|_{L^1_1(\mathbb{R}^3_v)}$ and
$\|g\|_{L\log L(\mathbb{R}^3_v)}$ such that in the case of $\gamma \in (0,1)$, there holds
\begin{align}
  \langle -Q(g,f),f \rangle \gtrsim C_g \|f\|^2_{H^s_{\frac{\gamma}{2}}}
                           - \left( \|g\|_{L^1_{2-\gamma}} + C_g \|g\|^2_{L^1} \right) \|f\|^2_{L^2_{\frac{\gamma}{2}}},
\end{align}
and in the case of $\gamma \in (-3,0]$, there holds
\begin{align}
  \langle -Q(g,f),f \rangle \gtrsim & C_g \|f\|^2_{H^s_{\frac{\gamma}{2}}}
                             - \left( \|g\|_{L^1_{|\gamma + 2|}} + \|g\|_{L^2_{|\gamma + 2|}} \right)
                                  \|f\|^2_{L^2_{\frac{\gamma}{2}}}
                             - \left( \|g\|_{L^1_{|\gamma|}} + \|g\|_{L^2_{|\gamma|}} \right)
                                  \|f\|^2_{H^\rho_{\frac{\gamma}{2}}},
\end{align}
with $-\frac{\gamma}{2} - \frac{3}{4} < \rho < s$.
\end{lemm}

{\noindent\bf Proof.} It is easy to check that
\begin{align}
  \langle Q(g,f),f \rangle = & \iiint_{\mathbb{R}^3 \times \mathbb{R}^3 \times \mathbb{S}^2} B g_* f(f'-f) dvdv_*d\sigma\\
 =&\frac{1}{2} \iiint_{\mathbb{R}^3 \times \mathbb{R}^3 \times \mathbb{S}^2} B g_* (f'^2-f^2) dvdv_*d\sigma
  -\frac{1}{2} \iiint_{\mathbb{R}^3 \times \mathbb{R}^3 \times \mathbb{S}^2} B g_* (f'-f)^2 dvdv_*d\sigma \nonumber\\
 \triangleq
  &\frac{1}{2} (\mathcal{J}_1-\mathcal{J}_2). \nonumber
\end{align}

\subsection{Upper bound for $\mathcal{J}_1$ }
By change of variables and the so-called cancellation lemma (see \cite{Entropy}), we can rewrite $\mathcal{J}_1 $ as
\begin{align}
\mathcal{J}_1 =& |\mathbb{S}^1| \iiint_{\mathbb{R}^3 \times \mathbb{R}^3 \times [0,\frac{\pi}{2}]}
                 sin \theta
                 \left\{ \frac{1}{cos^3\frac{\theta}{2}} B \left(\frac{|v-v_*|}{cos{\theta}/{2}},cos\theta \right)
                 -B(|v-v_*|,cos\theta)\right\} g_* f^2 dvdv_*d\theta \\
              =& |\mathbb{S}^1| \iiint_{\mathbb{R}^3 \times \mathbb{R}^3 \times [0,\frac{\pi}{2}]}
                 sin \theta b(cos\theta) |v-v_*|^\gamma
                 \left\{ \left(cos^{-1}\frac{\theta}{2}\right)^{\gamma+3}-1\right\}
                 g_* f^2 dvdv_*d\theta  \nonumber
\end{align}
Noticing that
$$
x^m-y^m=m\int_y^x z^{m-1} dz,
$$
we have, for $\theta \in [0,\frac{\pi}{2}]$,
\begin{align*}
  &\left(cos^{-1}\frac{\theta}{2}\right)^{\gamma+3}-1
  = \frac{\gamma+3}{cos^{\gamma+3}\frac{\theta}{2}}
    \int_0^{\frac{\theta}{2}} cos^{\gamma+2}\theta sin\theta d\theta  \\
  \leq &(\gamma+3)\int_0^{\frac{\theta}{2}} sin\theta d\theta
    \leq (\gamma+3)\left( 1-cos\frac{\theta}{2}\right) \nonumber\\
  \leq &2(\gamma+3) sin^2 \frac{\theta}{4} \leq C\theta^2 .\nonumber
\end{align*}
Then we arrive at
\begin{align}
  \mathcal{J}_1 \lesssim \iint_{\mathbb{R}^3 \times \mathbb{R}^3} |v-v_*|^\gamma g_* f^2 dvdv_*.
\end{align}

If $\gamma > 0$, we have
\begin{align} \label{J_1,0+}
  \mathcal{J}_1 \lesssim \|g\|_{L^1_{|\gamma|}} \|f\|^2_{L^2_{\frac{\gamma}{2}}}.
\end{align}

On the other aspect, if $\gamma \leq 0$, noticing the fact
\begin{align*}
  \left( \frac{|v-v_*|}{\langle v-v_* \rangle} \right)^{\gamma}
 \lesssim  {\sf 1}_{|v-v_*| \geq 1} + |v-v_*|^{\gamma} \cdot {\sf 1}_{|v-v_*| \leq 1},
\end{align*}
and
\begin{align*}
  \langle v-v_* \rangle^{\gamma} \lesssim \langle v_* \rangle^{-\gamma} \langle v \rangle^{\gamma},
\end{align*}
we can split $\mathcal{J}_1$ as
\begin{align}
  \mathcal{J}_1 \lesssim & \iint_{|v-v_*| \geq 1} g_* \langle v_* \rangle^{-\gamma}
                                           f^2 \langle v \rangle^{\gamma} dvdv_*
                   + \iint_{|v-v_*| \leq 1} |v-v_*|^\gamma g_* \langle v_* \rangle^{-\gamma}
                                             f^2 \langle v \rangle^{\gamma} dvdv_* \\
              \triangleq &  \mathcal{J}_{11} + \mathcal{J}_{12}.  \nonumber
\end{align}

Firstly we have
$$
 \mathcal{J}_{11} \lesssim \|g\|_{L^1_{|\gamma|}} \|f\|^2_{L^2_\frac{\gamma}{2}}.
$$

Furthermore, concerning with the estimate on $\mathcal{J}_{12}$, we cite the following result (see Lemma 2.6 in \cite{FiveGroup-III}):
\begin{lemm}\label{Lemma 2.6}
  If $\lambda < 3/2$ then
  \begin{align}
    \iint_{|v-v_*| \leq 1}  |f_*| \frac{|g|^2}{|v-v_*|^\lambda} dvdv_*
    \lesssim \|f\|_{L^2} \|g\|^2_{L^2}.
  \end{align}
  If $3/2 < \lambda < 3$ then
  \begin{align}
    \iint_{|v-v_*| \leq 1}  |f_*| \frac{|g|^2}{|v-v_*|^\lambda} dvdv_*
    \lesssim \|f\|_{L^2} \|g\|^2_{H^{\frac{\lambda}{2}-\frac{3}{4}}}.
  \end{align}
\end{lemm}

From the above lemma it follows that, if $\gamma > -3/2$,
$$
 \mathcal{J}_{12} \lesssim \|g\|_{L^2_{|\gamma|}} \|f\|^2_{L^2_{\frac{\gamma}{2}}},
$$
and if $-3 < \gamma \leq -3/2$,
$$
 \mathcal{J}_{12} \lesssim \|g\|_{L^2_{|\gamma|}} \|f\|^2_{H^\rho_{\frac{\gamma}{2}}},
$$
with $\rho > -\frac{\gamma}{2}-\frac{3}{4}$. Note that the case
$\gamma = -3/2$ can be treated as $\gamma -\varepsilon $ for any
small $\varepsilon > 0 $.

Thus we have, if $-3/2 < \gamma \leq 0$,
\begin{align}
  \mathcal{J}_1 \lesssim \left( \|g\|_{L^1_{|\gamma|}} + \|g\|_{L^2_{|\gamma|}} \right)
                         \|f\|^2_{L^2_{\frac{\gamma}{2}}},
\end{align}
and if  $-3 < \gamma \leq -3/2$,
\begin{align}
  \mathcal{J}_1 \lesssim \left( \|g\|_{L^1_{|\gamma|}} + \|g\|_{L^2_{|\gamma|}} \right)
                         \|f\|^2_{H^\rho_{\frac{\gamma}{2}}},
\end{align}
hence, together with the above two inequalities, we get for $\gamma \leq 0 $,
\begin{align}\label{J_1,-0}
  \mathcal{J}_1 \lesssim \left( \|g\|_{L^1_{|\gamma|}} + \|g\|_{L^2_{|\gamma|}} \right)
                         \|f\|^2_{H^\rho_{\frac{\gamma}{2}}},
\end{align}
where $\rho$ take the same value mentioned before.

\subsection{Lower bound for $\mathcal{J}_2$}

Setting $\mathcal{F}=f \langle v \rangle^{\frac{\gamma}{2}}$, we have
\begin{align}
  \mathcal{J}_2 = & \iiint_{\mathbb{R}^3 \times \mathbb{R}^3} |v-v_*|^\gamma b(cos\theta) g_* (f'-f)^2 dvdv_*d\sigma \\
                = & \iiint_{\mathbb{R}^3 \times \mathbb{R}^3} |v-v_*|^\gamma b(cos\theta) g_*
                     \left(\langle v' \rangle^{-\frac{\gamma}{2}}\mathcal{F}'-\langle v \rangle^{-\frac{\gamma}{2}}\mathcal{F} \right)^2 dvdv_*d\sigma. \nonumber
\end{align}
We shall give a different control for $\mathcal{J}_2$ with respect to $\gamma$.

\noindent $\bullet$ In the case of $\gamma \leq 0$:

Owing to $(A-B)^2 \geq \frac{A^2}{2}-B^2$, we obtain
\begin{align*}
  \left( \langle v' \rangle^{-\frac{\gamma}{2}}\mathcal{F}'-\langle v \rangle^{-\frac{\gamma}{2}}\mathcal{F} \right)^2
  \geq \frac{1}{2} \langle v' \rangle^{-\gamma}(\mathcal{F}'-\mathcal{F})^2
       -\mathcal{F}^2 \left( \langle v' \rangle^{-\frac{\gamma}{2}}-\langle v \rangle^{-\frac{\gamma}{2}} \right)^2,
\end{align*}
so we can write that
\begin{align}
  \mathcal{J}_2 \geq & \frac{1}{2} \iiint_{\mathbb{R}^3 \times \mathbb{R}^3} |v-v_*|^\gamma b(cos\theta) g_*
                                       \langle v' \rangle^{-\gamma}(\mathcal{F}'-\mathcal{F})^2 dvdv_*d\sigma \\
                     & -\iiint_{\mathbb{R}^3 \times \mathbb{R}^3} |v-v_*|^\gamma b(cos\theta) g_*
                              \mathcal{F}^2 \left( \langle v' \rangle^{-\frac{\gamma}{2}}
                                                  -\langle v \rangle^{-\frac{\gamma}{2}} \right)^2 dvdv_*d\sigma
                     \nonumber\\
               \triangleq & \mathcal{L}_1-\mathcal{L}_2. \nonumber
\end{align}

Observing that $|v-v_*| \thicksim |v'-v_*|$ and
$\langle v'-v_* \rangle^{-\gamma} \lesssim \langle v' \rangle^{-\gamma} \langle v_* \rangle^{-\gamma}$,
we have
$$
|v-v_*|^\gamma \langle v' \rangle^{-\gamma} \langle v_* \rangle^{-\gamma}
\gtrsim \left( \frac{|v-v_*|}{\langle v'-v_* \rangle} \right)^\gamma \gtrsim 1,
$$
then we get
\begin{align}
  \mathcal{L}_1 \gtrsim \iiint_{\mathbb{R}^3 \times \mathbb{R}^3 \times \mathbb{S}^2} b(cos\theta)
                              \left( g_* \langle v_* \rangle^\gamma \right)
                              \left( \mathcal{F}'-\mathcal{F} \right)^2 dvdv_*d\sigma.
\end{align}

Due to the well-known entropy dissipation inequality, we obtain
\begin{align}
  \mathcal{L}_1 \gtrsim C_g \| \mathcal{F} \|^2_{H^s},
\end{align}
where $C_g$ depends only on $\| g \langle v \rangle^\gamma \|_{L^1_1},~\| g \langle v \rangle^\gamma \|_{LlogL},$ and $b$.

As for the estimate of $\mathcal{L}_2$, we use the Taylor expansion at order $1$ to get:
\begin{align*}
   \langle v' \rangle^{-\frac{\gamma}{2}}-\langle v \rangle^{-\frac{\gamma}{2}}
  =\int_0^1 (v'-v) \nabla \left( \langle v \rangle^{-\frac{\gamma}{2}} \right) (v_\tau) d\tau,
\end{align*}
with $\tau \in [0,1]$ and $v_\tau=v+\tau(v'-v)$.

Since
\begin{align*}
   \left|(v'-v) \left( \nabla \big( \langle v \rangle^{-\frac{\gamma}{2}} \big) \right) (v_\tau) \right|
   \lesssim |v'-v| \langle v_\tau \rangle^{-\frac{\gamma}{2}-1}
   \lesssim \theta |v-v_*| \langle v_\tau \rangle^{-\frac{\gamma}{2}-1},
\end{align*}
then we have
\begin{align}
  \mathcal{L}_2 \lesssim & \iiint_{\mathbb{R}^3 \times \mathbb{R}^3 \times \mathbb{S}^2}
                              |v-v_*|^{\gamma+2} b(cos\theta) \theta^2 g_* \mathcal{F}^2
                              \langle v_\tau \rangle^{(-\gamma-2)} dvdv_*d\sigma \\
                \lesssim & \iint_{\mathbb{R}^3 \times \mathbb{R}^3} |v-v_*|^{\gamma+2} g_* \mathcal{F}^2
                              \langle v_\tau \rangle^{-(\gamma+2)} dvdv_* \nonumber.
\end{align}

Similar as we discussed before, if $\gamma +2 > 0$, we get
\begin{align*}
  |v-v_*|^{\gamma+2} \langle v_\tau \rangle^{-(\gamma+2)} \langle v_* \rangle^{-(\gamma+2)} \lesssim 1,
\end{align*}
which yields
\begin{align}
\mathcal{L}_2 \lesssim \iint_{\mathbb{R}^3 \times \mathbb{R}^3} g_* \mathcal{F}^2 \langle v_* \rangle^{\gamma+2}dvdv_*
\lesssim \| g \|_{L^1_{\gamma+2}} \|\mathcal{F}\|^2_{L^2}.
\end{align}

On the other hand, if $\gamma +2 \leq 0$, recalling the assumption
$\gamma \in (-3,0)$ immediately yields $\gamma +2 \in (-1,0]$,
then
\begin{align*}
  \left( \frac{|v-v_*|}{\langle v-v_* \rangle} \right)^{\gamma+2}
 \lesssim {\sf 1}_{|v-v_*| \geq 1} + |v-v_*|^{\gamma+2} \cdot {\sf 1}_{|v-v_*| \leq 1}.
\end{align*}

Noticing the fact
\begin{align*}
  \langle v-v_* \rangle^{\gamma+2} \langle v_\tau \rangle^{-(\gamma+2)} \langle v_* \rangle^{\gamma+2} \lesssim 1,
\end{align*}
then we obtain
\begin{align*}
  \mathcal{L}_2 \lesssim &\iint_{|v-v_*| \geq 1} g_* \mathcal{F}^2 \langle v_* \rangle^{-(\gamma+2)}dvdv_*
   +\iint_{|v-v_*| \leq 1} |v-v_*|^{\gamma+2} g_* \mathcal{F}^2 \langle v_* \rangle^{-(\gamma+2)}dvdv_*\\
  \triangleq & \mathcal{L}_{21} + \mathcal{L}_{22}. \nonumber
\end{align*}

Obviously we have
\begin{align*}
  \mathcal{L}_{21} \lesssim \|g\|_{L^1_{|\gamma +2|}} \|\mathcal{F}\|^2_{L^2}.
\end{align*}

Using Lemma \ref{Lemma 2.6} with $\gamma+2>-1>-3/2$, we get
\begin{align}
  \mathcal{L}_{22} \lesssim \|g\|_{L^2_{|\gamma +2|}} \|\mathcal{F}\|^2_{L^2},
\end{align}
so we arrive at the estimate
\begin{align}
  \mathcal{L}_2 \lesssim \left( \|g\|_{L^1_{|\gamma +2|}} + \|g\|_{L^2_{|\gamma +2|}} \right) \|\mathcal{F}\|^2_{L^2}.
\end{align}

Thus we have for $\gamma \leq 0$,
\begin{align}\label{J_2,-0}
  \mathcal{J}_2 = \mathcal{L}_1 - \mathcal{L}_2
               \gtrsim C_g \| \mathcal{F} \|^2_{H^s} - \left( \|g\|_{L^1_{|\gamma +2|}} + \|g\|_{L^2_{|\gamma +2|}}
                                                         \right) \|\mathcal{F}\|^2_{L^2},
\end{align}
where $C_g$ is defined before.

{\noindent $\bullet$ In the case of $\gamma > 0$:}

We first mention that
\begin{align*}
  \langle v \rangle \thicksim |v| + {\sf 1}_{|v| \leq 1},
\end{align*}
so
\begin{align}
  \mathcal{J}_2 \gtrsim \iiint_{\mathbb{R}^3 \times \mathbb{R}^3 \times \mathbb{S}^2}
                              \langle v-v_* \rangle^\gamma b(cos\theta) g_* (f'-f)^2 dvdv_*d\sigma
                        -\iiint_{\mathbb{R}^3 \times \mathbb{R}^3 \times \mathbb{S}^2}
                              b(cos\theta) g_* (f'-f)^2 dvdv_*d\sigma .
\end{align}

The fact $\langle v-v_* \rangle^\gamma \geq \langle v_* \rangle^{-\gamma} \langle v \rangle^\gamma $ implies
\begin{align*}
  (f'-f)^2 =  \left(\langle v' \rangle^{-\frac{\gamma}{2}}\mathcal{F}'-\langle v \rangle^{-\frac{\gamma}{2}}\mathcal{F}
                \right)^2
        \geq  \frac{1}{2} \langle v' \rangle^{-\gamma} (\mathcal{F}'-\mathcal{F})^2
            - \mathcal{F}^2 \left( \langle v' \rangle^{-\frac{\gamma}{2}}
                                  - \langle v \rangle^{-\frac{\gamma}{2}}\right)^2,
\end{align*}
from which we get
\begin{align}
  \mathcal{J}_2 \gtrsim & \frac{1}{2} \iiint_{\mathbb{R}^3 \times \mathbb{R}^3 \times \mathbb{S}^2}
                              b(cos\theta) \left( g_* \langle v_* \rangle^{-\gamma} \right)
                              (\mathcal{F}'-\mathcal{F})^2 dvdv_*d\sigma \nonumber \\
                        &- \iiint_{\mathbb{R}^3 \times \mathbb{R}^3 \times \mathbb{S}^2}
                              b(cos\theta)\left( g_* \langle v_* \rangle^{-\gamma} \right)
                              f^2 \left( \langle v' \rangle^{\frac{\gamma}{2}} - \langle v \rangle^{\frac{\gamma}{2}} \right)^2 dvdv_*d\sigma \nonumber \\
                        &- \iiint_{\mathbb{R}^3 \times \mathbb{R}^3 \times \mathbb{S}^2}
                              {\LARGE \bf 1}_{|v-v_*| \leq 1} b(cos\theta) g_* (f'-f)^2 dvdv_*d\sigma \nonumber \\
                \triangleq & \mathcal{L}_3 - \mathcal{L}_4 - \mathcal{L}_5. \nonumber
\end{align}

Similar as $\mathcal{L}_1$, we obtain
\begin{align}
  \mathcal{L}_3 \gtrsim C_g \| \mathcal{F}\|^2_{H^s},
\end{align}
with $C_g$ depends on $\| g \langle v \rangle^{-\gamma} \|_{L^1_1},~\| g \langle v \rangle^{-\gamma} \|_{LlogL},$ and $b$.

To bound $\mathcal{L}_4 $ we use the Taylor expansion for $ \langle v' \rangle^{\frac{\gamma}{2}} - \langle v \rangle^{\frac{\gamma}{2}} $, and get
\begin{align}
  \mathcal{L}_4 \lesssim \iiint_{\mathbb{R}^3 \times \mathbb{R}^3 \times \mathbb{S}^2}
                               b(cos\theta) g_* f^2 |v-v'|^2 \langle v_\tau \rangle^{\gamma -2}
                               \langle v_* \rangle^{-\gamma} dvdv_*d\sigma,
\end{align}
where $ \tau \in [0,1]$ and $v_\tau$ is defined as before.

Thus we obtain
\begin{align}
  \mathcal{L}_4 \lesssim & \iiint_{\mathbb{R}^3 \times \mathbb{R}^3 \times \mathbb{S}^2}
                               b(cos\theta) \theta^2 |v-v_*|^2 g_* f^2  \langle v_\tau \rangle^{\gamma -2}
                               \langle v_* \rangle^{-\gamma} dvdv_*d\sigma \\
                \lesssim & \iint_{\mathbb{R}^3 \times \mathbb{R}^3 } \left( g_* \langle v_* \rangle^{|\gamma-2|} \right)
                               \left( f^2 \langle v \rangle^{\gamma} \right) |v-v_*|^2
                                \langle v_* \rangle^{-\gamma} \langle v_* \rangle^{-|\gamma-2|}
                               \langle v \rangle^{-\gamma} \langle v_\tau \rangle^{\gamma -2} dvdv_* \nonumber \\
                \lesssim & \iint_{\mathbb{R}^3 \times \mathbb{R}^3 } \left( g_* \langle v_* \rangle^{|\gamma-2|} \right)
                               \left( f^2 \langle v \rangle^{\gamma} \right) |v-v_*|^2
                               \langle v-v_* \rangle^{-\gamma + \gamma-2} dvdv_*d\sigma \nonumber \\
                \lesssim & \|g\|_{L^1_{|\gamma -2 |}} \|\mathcal{F}\|^2_{L^2}. \nonumber
\end{align}
where we have used the following facts
\begin{align*}
  & \langle v-v_* \rangle^\gamma \lesssim \langle v_* \rangle^\gamma \langle v \rangle^\gamma,\\
  & \langle v_*-v_\tau \rangle^{-(\gamma-2)} \lesssim \langle v_* \rangle^{|\gamma-2|} \langle v_\tau\rangle^{-(\gamma-2)},
\end{align*}
and
$$|v_*-v_\tau| \thicksim |v-v_*|. $$

For the estimate of $\mathcal{L}_5$, we write it as
\begin{align}
  \mathcal{L}_5 = & \iiint_{\mathbb{R}^3 \times \mathbb{R}^3 \times \mathbb{S}^2}
                              {\Large \bf 1}_{|v-v_*| \leq 1} b(cos\theta) g_* (f'-f)^2 dvdv_*d\sigma\\
  = & -2 \iiint_{\mathbb{R}^3 \times \mathbb{R}^3 \times \mathbb{S}^2}
                              {\Large \bf 1}_{|v-v_*| \leq 1} b(cos\theta) g_* f(f'-f) dvdv_*d\sigma \nonumber\\
   & + \iiint_{\mathbb{R}^3 \times \mathbb{R}^3 \times \mathbb{S}^2}
                              {\Large \bf 1}_{|v-v_*| \leq 1} b(cos\theta) g_* (f'^2-f^2) dvdv_*d\sigma\nonumber\\
  \triangleq & \mathcal{L}^{(1)}_5 + \mathcal{L}^{(2)}_5. \nonumber
\end{align}

When bounding the term $\mathcal{L}^{(1)}_5$ we refer to Proposition 2.1 and Remark 2.2 in \cite{FiveGroup-III}, then we have
\begin{align}
  |\mathcal{L}^{(1)}_5| \lesssim \|g\|_{L^1} \|f\|^2_{H^s}.
\end{align}

Under a similar argument as we did for $\mathcal{J}_1$ we get
\begin{align}
  \mathcal{L}^{(2)}_5 \lesssim \|g\|_{L^1} \|f\|^2_{L^2},
\end{align}
so we have the following estimate
\begin{align}
  \mathcal{L}_5 \lesssim & \|g\|_{L^1} \left( \|f\|^2_{L^2} + \|f\|^2_{H^s} \right) \\
                \lesssim & \|g\|_{L^1} \left( \|f\|^2_{L^2} + \|f\|^2_{\dot{H}^s} \right).\nonumber
\end{align}

Therefore we obtain for $\gamma > 0$,
\begin{align}
  \mathcal{J}_2 \gtrsim & \mathcal{L}_3 - \mathcal{L}_4 - \mathcal{L}_5 \\
                \gtrsim & C_g \|\mathcal{F}\|^2_{H^s} - \|g\|_{L^1_{|\gamma -2 |}} \|\mathcal{F}\|^2_{L^2}
                          - \|g\|_{L^1} \big( \|f\|^2_{L^2} + \|f\|^2_{\dot{H}^s} \big) \nonumber \\
                \gtrsim & C_g \|\mathcal{F}\|^2_{H^s} - \|g\|_{L^1_{|\gamma -2 |}} \|\mathcal{F}\|^2_{L^2}
                          - \|g\|_{L^1} \|f\|^2_{\dot{H}^s}. \nonumber
\end{align}

Following along the same lines as in \cite{Chen-He-1}, we finally arrive at
\begin{align}\label{J_2,0+}
  \mathcal{J}_2 \gtrsim  C_g \|\mathcal{F}\|^2_{H^s}
                        - \left( \|g\|_{L^1_{|\gamma -2 |}} + C_g \|g\|^2_{L^1} \right)
                            \|\mathcal{F}\|^2_{L^2}.
\end{align}

Now we conclude that (\ref{J_1,0+}) and (\ref{J_2,0+}) imply that for $\gamma \in (0,1)$,
\begin{align}
  \langle -Q(g,f),f \rangle \gtrsim C_g \|f\|^2_{H^s_{\frac{\gamma}{2}}}
                           - \left( \|g\|_{L^1_{2-\gamma}} + C_g \|g\|^2_{L^1} \right) \|f\|^2_{L^2_{\frac{\gamma}{2}}},
\end{align}
where we have used $|\gamma - 2|=2-\gamma >\gamma$.

On the other hand, (\ref{J_1,-0}) and (\ref{J_2,-0}) imply that
for $\gamma \leq 0$,
\begin{align}
  \langle -Q(g,f),f \rangle \gtrsim & C_g \|f\|^2_{H^s_{\frac{\gamma}{2}}}
                             - \left( \|g\|_{L^1_{|\gamma + 2|}} + \|g\|_{L^2_{|\gamma + 2|}} \right)
                                  \|f\|^2_{L^2_{\frac{\gamma}{2}}}
                             - \left( \|g\|_{L^1_{|\gamma|}} + \|g\|_{L^2_{|\gamma|}} \right)
                                  \|f\|^2_{H^\rho_{\frac{\gamma}{2}}},
\end{align}
with $-\frac{\gamma}{2} - \frac{3}{4} < \rho < s$.

This leads us to the conclusion. \qed

\section{Proof of the main lemma}
\subsection{Rewrite the equation}
Let $\mu =\mu_{\delta,\kappa}(t) =e^{-(\delta-\kappa t) \langle v \rangle ^2}  $ with $\delta > \kappa T$. According to the translation invariance of the collision operator with respect to the variable $v$ (see \cite{Desvillettes2, Ukai-84}), we have, for the translation operator $\tau_h $ in $v$ by $h$,
$$
\tau_h Q(f,g)=Q(\tau_h f,\tau_h g).
$$

Thus we have
$$
\partial_v^\alpha Q(f,g) = \sum_{\alpha=\alpha'+\alpha''}
                           \frac{\alpha!}{\alpha'! \alpha''!} Q\big(f^{(\alpha')},g^{(\alpha'')}\big).
$$

Then we could obtain from Eq.(\ref{BE}) that
$$
\partial_t (\partial_v^\alpha f) = Q(f,f^{(\alpha)})+\sum_{\alpha' \neq 0}
                            \frac{\alpha!}{\alpha'! \alpha''!}Q\big(f^{(\alpha')},f^{(\alpha'')}\big).
$$

Multiplying both sides by $\mu^{-1}$, we obtain
\begin{align}\label{eq1}
(\partial_t +\kappa \langle v \rangle ^2)(\mu^{-1} \partial_v^\alpha f)= \mu^{-1}Q(f,f^{(\alpha)})
 +\sum_{\alpha' \neq 0}\frac{\alpha!}{\alpha'! \alpha''!}\mu^{-1} Q\big(f^{(\alpha')},f^{(\alpha'')}\big).
\end{align}

Set $F=\mu^{-1}f$ and denote $F^{(\alpha)}=\mu^{-1}f^{(\alpha)}$ for $\alpha \in \mathbb{Z}^3_+$. Noticing that $\mu \mu_* = \mu' \mu'_*$, we get the following formula
$$
\mu^{-1}Q(f,g)=Q(\mu F,G) + \iint B(\mu_*-\mu'_*)F'_* G'dv_* d\sigma.
$$

Then it follows from (\ref{eq1}) that
\begin{align*}
(\partial_t +\kappa \langle v \rangle ^2) F^{(\alpha)}=&Q(\mu F,F^{(\alpha)})
+\sum_{\alpha' \neq 0}\frac{\alpha!}{\alpha'! \alpha''!}Q\big(\mu F^{(\alpha')},F^{(\alpha'')}\big)\\
&+\sum_{\alpha=\alpha'+\alpha''}\frac{\alpha!}{\alpha'! \alpha''!}
                       \iint B(\mu_*-\mu'_*)(F^{(\alpha')})'_* (F^{(\alpha'')})'dv_* d\sigma.
\end{align*}

Hereafter we denote $W_l=\langle v \rangle ^l$. Multiplying by $W_l^2F^{(\alpha)}$ both sides and integrating over $v$, we have
\begin{align}\label{EQ}
&\frac{1}{2}\frac{d}{dt}\|W_l F^{(\alpha)}\|^2 + \kappa \|W_{l+1} F^{(\alpha)}\|^2\\
=&\big\langle  Q(\mu F,F^{(\alpha)}), W_l^2F^{(\alpha)}\big\rangle
   +\sum_{\alpha' \neq 0}\frac{\alpha!}{\alpha'!\alpha''!}
    \Big\langle Q\big(\mu F^{(\alpha')},F^{(\alpha'')}\big),W_l^2F^{(\alpha)}\Big\rangle \nonumber\\
&+\sum_{\alpha=\alpha'+\alpha''}\frac{\alpha!}{\alpha'! \alpha''!}
     \iiint B(\mu_*-\mu'_*)(F^{(\alpha')})'_* (F^{(\alpha'')})'W_l^2F^{(\alpha)}dv dv_* d\sigma \nonumber\\
=&\Psi_1^{(0,\alpha)}(t)
  +\sum_{\alpha' \neq 0}\frac{\alpha!}{\alpha'! \alpha''!}\Psi_1^{(\alpha',\alpha'')}(t)
  +\sum_{\alpha=\alpha'+\alpha''}\frac{\alpha!}{\alpha'! \alpha''!}\Psi_2^{(\alpha',\alpha'')}(t)
  \nonumber\\
=&\Psi_1^{(0,\alpha)}(t)+\mathcal{J}(t)+\mathcal{K}(t)\nonumber.
\end{align}

Then multiplying by the weight $\frac{\rho^{2|\alpha|}}{\{(\alpha-r)!\}^{2\nu}}$ both sides, and
integrating from $0$ to $t\in (0,T]$, we obtain
\begin{align}\label{EQ2}
&\|f(t)\|^2_{\delta -\kappa t,l,\rho,\alpha,r}
+ 2\kappa \int^t_0  \|f(\tau)\|^2_{\delta-\kappa \tau,l+1,\rho,\alpha,r}d\tau \\
\leq
&\|f(0)\|^2_{\delta,l,\rho,\alpha,r}
 +2\int^t_0 \frac{\rho^{2|\alpha|}}{\{(\alpha-r)!\}^{2\nu}}
   \big( \Psi_1^{(0,\alpha)}(\tau)+\mathcal{J}(\tau)+\mathcal{K}(\tau) \big)d\tau
   \nonumber.
\end{align}

We will estimate the above terms one by one. Firstly we consider the estimate on the ``remainder term"
$\Psi_2^{(\alpha',\alpha'')}(t)$, which will be given in the next subsection.

\subsection{The estimate on the remainder term}

We rewrite the remainder term $\Psi_2^{(\alpha',\alpha'')}(t)$ as follows:
\begin{align*}
&\Psi_2^{(\alpha',\alpha'')}(t)\\
=&\iiint B(\mu_*-\mu'_*)(F^{(\alpha')})'_* (W_{l-1}F^{(\alpha'')})'W_{l+1}F^{(\alpha)}dv dv_* d\sigma\\
&+\iiint B(\mu_*-\mu'_*)(F^{(\alpha')})'_* \big(W_{l-1}-W'_{l-1}\big)
  (F^{(\alpha'')})'W_{l+1}F^{(\alpha)}dv dv_* d\sigma\\
=&\Psi_{2,1}^{(\alpha',\alpha'')}(t)+\Psi_{2,2}^{(\alpha',\alpha'')}(t).
\end{align*}

\vspace*{1em}
{\noindent \it Step I: The estimate for $\Psi_{2,2}^{(\alpha',\alpha'')}(t)$. }
The term $\Psi_{2,2}^{(\alpha',\alpha'')}(t)$ can be written as
\begin{align}
  \Psi_{2,2}^{(\alpha',\alpha'')}(t)
= &\iiint B \mu_* (F^{(\alpha')})'_* \big(W_{l-1}-W'_{l-1}\big)
  (F^{(\alpha'')})'W_{l+1}F^{(\alpha)}dv dv_* d\sigma \\
 & - \iiint B \mu'_* (F^{(\alpha')})'_* \big(W_{l-1}-W'_{l-1}\big)
  (F^{(\alpha'')})'W_{l+1}F^{(\alpha)}dv dv_* d\sigma \nonumber \\
= &I_+ - I_-. \nonumber
\end{align}

The Taylor expansion up to order 2 gives
\begin{align*}
  W_{l-1}-W'_{l-1} = (v-v')\cdot (\nabla_v W_{l-1})(v') + \int_0^1 (1-\tau)(v-v')^2 (\nabla^2_v W_{l-1})(v_\tau) d\tau
\end{align*}
with $\tau \in [0,1]$ and $v_\tau = v'+\tau(v-v')$.

Corresponding to the two terms of the right-hand side, we rewrite $I_+$ (and $I_-$) as $I_{+,1} + I_{+,2}$ (and $I_{-,1} + I_{-,2}$, respectively).

Furthermore, we have
\begin{align*}
  \left| (v-v')^2 (\nabla^2_v W_{l-1})(v_\tau) \right|
\lesssim & |v-v'|^2 \langle v_\tau \rangle^{l-3}
\lesssim   |v-v'|^2 \left ( \langle v_\tau -v_* \rangle^{l-3}
                         + \langle v_* \rangle^{l-3} \right ) \\
\lesssim & \theta^2 |v-v_*|^2 \left ( \langle v' \rangle^{l-3}
                                 + \langle v_* \rangle^{l-3} \right )
\lesssim \theta^2 |v-v_*|^2 \langle v' \rangle^{l-3}
                                  \langle v_* \rangle^{l-3}, \\
\Big(or~ \lesssim  & \theta^2 |v-v_*|^2 \left ( \langle v' \rangle^{l-3}
                                 + \langle v'_* \rangle^{l-3} \right )
\lesssim \theta^2 |v-v_*|^2 \langle v' \rangle^{l-3}
                                  \langle v'_* \rangle^{l-3},\Big)
\end{align*}
where we have used the fact $ |v'-v_*| \thicksim |v-v_*| $ and $|v_*|^2 \leq |v|^2 + |v_*|^2 = |v'|^2 + |v'_*|^2. $

Observing the fact $|\mu_* W_{l-3,*}|, ~|\mu'_* W'_{l-3,*}| \lesssim 1$, we get
\begin{align}
  | I_{+,2} |
\lesssim & \iiint_{\mathbb{R}^3 \times \mathbb{R}^3 \times \mathbb{S}^2}
|v-v_*|^{\gamma+2} b(cos\theta) \theta^2 |\mu_* W_{l-3,*}| |\left( F^{(\alpha')} \right)'_*|
|\left( W_{l-3} F^{(\alpha'')} \right)'|
|\left( W_{l+1}F^{(\alpha)} \right)| dv dv_* d\sigma \\
\lesssim & \iiint_{\mathbb{R}^3 \times \mathbb{R}^3 \times \mathbb{S}^2}
|v-v_*|^{\gamma+2} b(cos\theta) \theta^2 |\left( F^{(\alpha')} \right)'_*|
|\left( W_{l-3} F^{(\alpha'')} \right)'| |\left( W_{l+1}F^{(\alpha)} \right)|
dv dv_* d\sigma \nonumber\\
\lesssim & \iiint_{\mathbb{R}^3 \times \mathbb{R}^3 \times \mathbb{S}^2}
|v-v_*|^{\gamma+2} b(cos\theta) \theta^2 |\left( F^{(\alpha')} \right)_*|
|\left( W_{l-3} F^{(\alpha'')} \right)| |\left( W_{l+1}F^{(\alpha)} \right)'|
dv dv_* d\sigma \nonumber\\
\triangleq &I_2,\nonumber
\end{align}
and similarly,
\begin{align}
  | I_{-,2} |
\lesssim & \iiint_{\mathbb{R}^3 \times \mathbb{R}^3 \times \mathbb{S}^2}
|v-v_*|^{\gamma+2} b(cos\theta) \theta^2 |\mu'_* W'_{l-3,*}| |\left( F^{(\alpha')} \right)'_*|
|\left( W_{l-3} F^{(\alpha'')} \right)'|
|\left( W_{l+1}F^{(\alpha)} \right)| dv dv_* d\sigma \\
\lesssim & \iiint_{\mathbb{R}^3 \times \mathbb{R}^3 \times \mathbb{S}^2}
|v-v_*|^{\gamma+2} b(cos\theta) \theta^2 |\left( F^{(\alpha')} \right)'_*|
|\left( W_{l-3} F^{(\alpha'')} \right)'| |\left( W_{l+3}F^{(\alpha)} \right)|
dv dv_* d\sigma \nonumber\\
\lesssim & \iiint_{\mathbb{R}^3 \times \mathbb{R}^3 \times \mathbb{S}^2}
|v-v_*|^{\gamma+2} b(cos\theta) \theta^2 |\left( F^{(\alpha')} \right)_*|
|\left( W_{l-3} F^{(\alpha'')} \right)| |\left( W_{l+1}F^{(\alpha)} \right)'|
dv dv_* d\sigma \nonumber\\
= & I_2. \nonumber
\end{align}

Now it suffices to bound the integral $I_2$. If $\gamma+2 \geq 0$, noticing that
\begin{align*}
  |v-v_*|^{\gamma+2} \lesssim \langle v_\tau \rangle^{\gamma+2} \langle v_* \rangle^{\gamma+2},
\end{align*}
we have
\begin{align}
   I_2
\lesssim & \iiint_{\mathbb{R}^3 \times \mathbb{R}^3 \times \mathbb{S}^2}
               |v-v_*|^{\gamma+2} b(cos\theta) \theta^2 |\left( F^{(\alpha')} \right)_*|
               |\left( W_{l-3} F^{(\alpha'')} \right)| |\left( W_{l+1}F^{(\alpha)} \right)'|
               dv dv_* d\sigma \\
\lesssim & \iiint_{\mathbb{R}^3 \times \mathbb{R}^3 \times \mathbb{S}^2} b(cos\theta) \theta^2
               |\left( W_{\gamma+2} F^{(\alpha')} \right)_*|
               |\left( W_{l-1+\gamma} F^{(\alpha'')} \right)| |\left( W_{l+1}F^{(\alpha)} \right)'|
               dv dv_* d\sigma \nonumber\\
\lesssim & \left( \iiint_{\mathbb{R}^3 \times \mathbb{R}^3 \times \mathbb{S}^2} b(cos\theta) \theta^2
               |\left( W_{\gamma+2} F^{(\alpha')} \right)_*|
               |\left( W_{l-1+\gamma} F^{(\alpha'')} \right)|^2 dv dv_* d\sigma \right)^{1/2} \nonumber\\
 & \times \left( \iiint_{\mathbb{R}^3 \times \mathbb{R}^3 \times \mathbb{S}^2} b(cos\theta) \theta^2
               |\left( W_{\gamma+2} F^{(\alpha')} \right)_*|
               |\left( W_{l+1}F^{(\alpha)} \right)'|^2 dv dv_* d\sigma \right)^{1/2} \nonumber\\
\triangleq & G_1^{1/2} \times G_2^{1/2}. \nonumber
\end{align}

The assumption $s \le 1/2$ implies that $\int_{\mathbb{S}^2} b(\cos \theta) \theta^2 d\sigma \leq C$. Thus
\begin{align*}
  G_1 \lesssim \|W_{\gamma+2}F^{(\alpha')} \|_{L^1} \|W_{l-1+\gamma}F^{(\alpha'')} \|^2_{L^2}.
\end{align*}

By virtue of the regular change of variables $v\mapsto v'=\frac{v+v_*}{2}+\frac{|v-v_*|}{2}\sigma$ for fixed $\sigma$ and $v_*$ whose Jacobian satisfying:
$$
\left|\frac{dv'}{dv}\right|=\frac{\cos^2(\theta/2)}{4},
$$
we have,
\begin{align*}
G_2 \lesssim & \iiint_{\mathbb{R}^3 \times \mathbb{R}^3 \times \mathbb{S}^2}
                     b(\cos \theta) \theta^2 \frac{4}{\cos^2(\theta/2)}
                     |\left( W_{\gamma+2} F^{(\alpha')} \right)_*|
                     |\left( W_{l+1}F^{(\alpha)} \right)'|^2 dv dv_* d\sigma\\
\lesssim & \|W_{\gamma+2} F^{(\alpha')} \|_{L^1} \|W_{l+1}F^{(\alpha)}\|^2_{L^2}.
\end{align*}

So, we obtain in the case of $\gamma+2 \geq 0$:
\begin{align}
  I_2 \lesssim & \|W_{\gamma+2} F^{(\alpha')} \|_{L^1} \|W_{l-1+\gamma} F^{(\alpha'')} \|_{L^2}
                 \|W_{l+1}F^{(\alpha)}\|_{L^2} \\
      \lesssim & \|W_l F^{(\alpha')} \|_{L^2} \|W_l F^{(\alpha'')} \|_{L^2}
                 \|W_{l+1}F^{(\alpha)}\|_{L^2}, \nonumber
\end{align}
for any large $l \geq 4$.

On the other hand, if $\gamma+2 < 0$, we have
\begin{align}
  I_2 \lesssim & \iiint_{|v-v_*| \geq 1} b(cos\theta) \theta^2
                     |\left( F^{(\alpha')} \right)_*|
                     |\left( W_{l-3} F^{(\alpha'')} \right)| |\left( W_{l+1}F^{(\alpha)} \right)'|
                     dv dv_* d\sigma \\
            & + \iiint_{|v-v_*| \leq 1} |v-v_*|^{\gamma+2} b(cos\theta) \theta^2
                     |\left( F^{(\alpha')} \right)_*|
                     |\left( W_{l-3} F^{(\alpha'')} \right)| |\left( W_{l+1}F^{(\alpha)} \right)'|
                     dv dv_* d\sigma \nonumber\\
  \lesssim & \left( \iiint_{|v-v_*| \geq 1} b(cos\theta) \theta^2
                     |\left( F^{(\alpha')} \right)_*| |\left( W_{l-3} F^{(\alpha'')} \right)|^2
                     dv dv_* d\sigma \right)^{1/2} \nonumber\\
  & ~ \times \left( \iiint_{|v-v_*| \geq 1} b(cos\theta) \theta^2
                     |\left( F^{(\alpha')} \right)_*| |\left( W_{l+1}F^{(\alpha)} \right)'|^2
                     dv dv_* d\sigma \right)^{1/2} \nonumber\\
  & +        \left( \iiint_{|v-v_*| \leq 1} |v-v_*|^{\gamma+2} b(cos\theta) \theta^2
                     |\left( F^{(\alpha')} \right)_*| |\left( W_{l-3} F^{(\alpha'')} \right)|^2
                     dv dv_* d\sigma \right)^{1/2} \nonumber\\
  & ~ \times \left( \iiint_{|v-v_*| \leq 1} |v-v_*|^{\gamma+2} b(cos\theta) \theta^2
                     |\left( F^{(\alpha')} \right)_*| |\left( W_{l+1}F^{(\alpha)} \right)'|^2
                     dv dv_* d\sigma \right)^{1/2} \nonumber\\
\triangleq & G_3^{1/2} \times G_4^{1/2} + G_5^{1/2} \times G_6^{1/2}.\nonumber
\end{align}

Firstly we have
\begin{align*}
  & G_3 \lesssim \| F^{(\alpha')} \|_{L^1} \|W_{l-3}F^{(\alpha'')}\|^2_{L^2},\\
  & G_4 \lesssim \| F^{(\alpha')} \|_{L^1} \|W_{l+1}F^{(\alpha)}\|^2_{L^2}.
\end{align*}

Furthermore, due to the fact $\gamma+2 \in (-1/2,0)$ and Lemma \ref{Lemma 2.6}, we get
\begin{align*}
  & G_5 \lesssim \| F^{(\alpha')} \|_{L^2} \|W_{l-3}F^{(\alpha'')}\|^2_{L^2},\\
  & G_6 \lesssim \| F^{(\alpha')} \|_{L^2} \|W_{l+1}F^{(\alpha)}\|^2_{L^2}.
\end{align*}

Combining the above estimates for $G_3,~G_4,~G_5$, and $G_6$, we obtain
\begin{align}
  I_2 \lesssim & \left( \| F^{(\alpha')} \|_{L^1} + \| F^{(\alpha')} \|_{L^2} \right)
                  \|W_{l-1+\gamma} F^{(\alpha'')}\|_{L^2}
                  \|W_{l+1}F^{(\alpha)}\|_{L^2} \\
      \lesssim & \|W_l F^{(\alpha')} \|_{L^2} \|W_l F^{(\alpha'')} \|_{L^2}
                 \|W_{l+1}F^{(\alpha)}\|_{L^2}, \nonumber
\end{align}
for any $l \geq 4$.

Concerning the estimate on $I_{+,1}$ (or $I_{-,1}$), we use the symmetry of cross-section $B$ with respect to $\sigma$ around the direction $(v-v_*)/|v-v_*|$ (see \cite{FiveGroup-regulariz, Gressman}), which forces all components of $v-v'$ to vanish except the component in the symmetry direction. Noticing $(v-v') \perp (v'-v_*)$, we can take place of $v-v'$ in $I_{+,1}$ (or $I_{-,1}$) by
\begin{align*}
    \left \langle v-v', \frac{v-v_*}{|v-v_*|} \right \rangle \cdot \frac{v-v_*}{|v-v_*|}
= & \left \langle v-v', \frac{v-v'+v'-v_*}{|v-v_*|} \right \rangle \cdot \frac{v-v_*}{|v-v_*|} \\
= &(v-v_*) \frac{|v-v'|^2}{|v-v_*|^2} \\
= &(v-v_*) sin^2 \frac{\theta}{2}.
\end{align*}

Since
\begin{align*}
  \left| (v-v_*)\cdot (\nabla_v W_{l-1})(v') \right| \lesssim |v-v_*| \langle v' \rangle^{l-2},
\end{align*}
then we have
\begin{align}
  |I_{+,1}|
  \lesssim & \iiint_{\mathbb{R}^3 \times \mathbb{R}^3 \times \mathbb{S}^2}
               b(cos \theta) \theta^2 |v-v_*|^{\gamma+1} |\left( F^{(\alpha')} \right)'_*|
               |\left( W_{l-2} F^{(\alpha'')} \right)'| |\left( W_{l+1}F^{(\alpha)} \right)|
               dv dv_* d\sigma \\
  \lesssim & \iiint_{\mathbb{R}^3 \times \mathbb{R}^3 \times \mathbb{S}^2}
               b(cos \theta) \theta^2 |v-v_*|^{\gamma+1} |\left( F^{(\alpha')} \right)_*|
               |\left( W_{l-2} F^{(\alpha'')} \right)| |\left( W_{l+1}F^{(\alpha)} \right)'|
               dv dv_* d\sigma \nonumber\\
  \triangleq & I_1, \nonumber
\end{align}

Arguing as we treated with $I_2$, we have, if $\gamma+1 \geq 0$,
\begin{align}
  I_1 \lesssim & \|W_{\gamma+1} F^{(\alpha')} \|_{L^1} \|W_{l-1+\gamma} F^{(\alpha'')} \|_{L^2}
                 \|W_{l+1}F^{(\alpha)}\|_{L^2} \\
      \lesssim & \|W_l F^{(\alpha')} \|_{L^2} \|W_l F^{(\alpha'')} \|_{L^2}
                 \|W_{l+1}F^{(\alpha)}\|_{L^2}, \nonumber
\end{align}
for any large $l \geq 4$. And if $\gamma+1 < 0$, using Lemma \ref{Lemma 2.6} with $\gamma+1 \in (-3/2,0)$, we have
\begin{align}
  I_1 \lesssim & \left( \|F^{(\alpha')} \|_{L^1} + \|F^{(\alpha')} \|_{L^2} \right)
                    \|W_{l-2} F^{(\alpha'')} \|_{L^2}
                    \|W_{l+1}F^{(\alpha)}\|_{L^2} \\
      \lesssim & \|W_l F^{(\alpha')} \|_{L^2} \|W_l F^{(\alpha'')} \|_{L^2}
                 \|W_{l+1}F^{(\alpha)}\|_{L^2}, \nonumber
\end{align}
for any $l \geq 4$.

Since we can deduce the result $I_{-,1} \lesssim I_1$, then we have for $\gamma \in (-5/2,1)$,
\begin{align}
  \left| \Psi_{2,2}^{(\alpha',\alpha'')} \right|
\leq & |I_{+,1}|+|I_{+,2}|+|I_{-,1}|+|I_{-,2}|
\lesssim I_1 + I_2 \\
\lesssim &\|W_l F^{(\alpha')} \|_{L^2} \|W_l F^{(\alpha'')} \|_{L^2}
                 \|W_{l+1}F^{(\alpha)}\|_{L^2}, \nonumber
\end{align}
for any $l \geq 4$.

\vspace*{2em}
{\noindent\it Step II: The estimate for $\Psi_{2,1}^{(\alpha',\alpha'')}(t)$. }
Thanks to the Taylor expansion
\begin{align*}
 \mu_*-\mu'_* = \int_0^1 (v_*-v'_*) \cdot \left( \nabla_v \mu \right)(v_\tau) d\tau,
\end{align*}
with $\tau \in [0,1]$ and $v_\tau = v'_* +\tau (v_*-v'_*)$. Since $|\nabla_v \mu| \lesssim 1$, we have
\begin{align}
  \Psi_{2,1}^{(\alpha',\alpha'')}(t)
\lesssim & \iiint_{\mathbb{R}^3 \times \mathbb{R}^3 \times \mathbb{S}^2}
           |v-v_*|^{\gamma+1}b(\cos\theta) \theta^2 |\left( F^{(\alpha')} \right)'_*|
           |\left( W_{l-1} F^{(\alpha'')} \right)'|
           |\left( W_{l+1} F^{(\alpha)} \right)| dvdv_*d\sigma \\
  =  & \iiint_{\mathbb{R}^3 \times \mathbb{R}^3 \times \mathbb{S}^2}
           |v-v_*|^{\gamma+1} b(\cos\theta) \theta^2 |\left( F^{(\alpha')} \right)_*|
           |\left( W_{l-1} F^{(\alpha'')} \right)|
           |\left( W_{l+1} F^{(\alpha)} \right)'| dvdv_*d\sigma, \nonumber
\end{align}
where by symmetry we have replaced $v_*-v'_*$ with
\begin{align*}
  \left \langle v_*-v'_*, \frac{v'-v'_*}{|v'-v'_*|} \right \rangle \cdot \frac{v'-v'_*}{|v'-v'_*|}
= (v'-v'_*) \cdot \frac{| v_*-v'_*|^2}{|v'-v'_*|^2}
= (v'-v'_*) sin^2 \frac{\theta}{2}.
\end{align*}

Thus, if $\gamma+1 \geq 0$, we have
\begin{align}
  \Psi_{2,1}^{(\alpha',\alpha'')}(t)
\lesssim & \iiint_{\mathbb{R}^3 \times \mathbb{R}^3 \times \mathbb{S}^2}
           b(\cos\theta) \theta^2 |W_{\gamma+1}\left( F^{(\alpha')} \right)'_*|
           |\left( W_{l+\gamma} F^{(\alpha'')} \right)'|
           |\left( W_{l+1} F^{(\alpha)} \right)| dvdv_*d\sigma \\
\leq & \left\{ \iiint_{\mathbb{R}^3 \times \mathbb{R}^3 \times \mathbb{S}^2}
           b(\cos\theta) \theta^2 |\left( W_{\gamma+1} F^{(\alpha')} \right)_*|
           |\left( W_{l+\gamma} F^{(\alpha'')} \right)|^2 dvdv_*d\sigma
       \right\}^{1/2} \nonumber \\
& \times \left\{ \iiint_{\mathbb{R}^3 \times \mathbb{R}^3 \times \mathbb{S}^2}
           b(\cos\theta) \theta^2  |\left( W_{\gamma+1}F^{(\alpha')} \right)_*|
           |\left( W_{l+1} F^{(\alpha)} \right)'|^2 dvdv_*d\sigma
         \right\}^{1/2} \nonumber \\
\triangleq & G_7^{1/2} \times G_8^{1/2}, \nonumber
\end{align}

We immediately have
\begin{align}
  G_7 \lesssim \|W_{\gamma+1}F^{(\alpha')} \|_{L^1} \|W_{l+\gamma}F^{(\alpha'')} \|^2_{L^2}.
\end{align}

Since the H\"{o}lder inequality yields
$$
\| W_{\eta} G\|_{L^2}\leq \|G\|^{1-\eta}_{L^2} \|W_1 G\|^\eta_{L^2},
$$
for $\eta \in [0,1]$, then, denoting $\eta \triangleq 1-\gamma^+$, we get
\begin{align}
  G_7 \lesssim \|W_{\gamma+1}F^{(\alpha')} \|_{L^1} \|W_{l}F^{(\alpha'')} \|^{2\eta}_{L^2}
                                           \|W_{l+1}F^{(\alpha'')} \|^{2(1-\eta)}_{L^2}.
\end{align}

As for the term $G_8$, by taking advantage of a regular change of variables $v \mapsto v'$, we get
\begin{align}
  G_8 \lesssim \|W_{\gamma+1}F^{(\alpha')} \|_{L^1} \|W_{l+1}F^{(\alpha)} \|^2_{L^2},
\end{align}
so we deduce that, for $\gamma+1 \geq 0$,
\begin{align}
  \Psi_{2,1}^{(\alpha',\alpha'')}(t)
\lesssim & \|W_{\gamma+1}F^{(\alpha')} \|_{L^1}  \|W_{l}F^{(\alpha'')} \|^{\eta}_{L^2}
           \|W_{l+1}F^{(\alpha'')} \|^{1-\eta}_{L^2} \|W_{l+1}F^{(\alpha)} \|_{L^2} \\
\lesssim & \|W_l F^{(\alpha')} \|_{L^2} \|W_{l}F^{(\alpha'')} \|^{\eta}_{L^2}
           \|W_{l+1}F^{(\alpha'')} \|^{1-\eta}_{L^2} \|W_{l+1}F^{(\alpha)} \|_{L^2}. \nonumber
\end{align}

If $\gamma+1 < 0$, we can write that
\begin{align}
  \Psi_{2,1}^{(\alpha',\alpha'')}(t)
\lesssim & \iiint_{|v-v_*| \geq 1} b(cos\theta) \theta^2
                     |\left( F^{(\alpha')} \right)_*|
                     |\left( W_{l-1} F^{(\alpha'')} \right)| |\left( W_{l+1}F^{(\alpha)} \right)'|
                     dv dv_* d\sigma \\
            & + \iiint_{|v-v_*| \leq 1} |v-v_*|^{\gamma+1} b(cos\theta) \theta^2
                     |\left( F^{(\alpha')} \right)_*|
                     |\left( W_{l-1} F^{(\alpha'')} \right)| |\left( W_{l+1}F^{(\alpha)} \right)'|
                     dv dv_* d\sigma \nonumber\\
\lesssim & \left\{ \iiint_{|v-v_*| \geq 1} b(cos\theta) \theta^2
                     |\left( F^{(\alpha')} \right)_*| |\left( W_{l-1} F^{(\alpha'')} \right)|^2
                     dv dv_* d\sigma \right\}^{1/2} \nonumber\\
  & ~ \times \left\{ \iiint_{|v-v_*| \geq 1} b(cos\theta) \theta^2
                     |\left( F^{(\alpha')} \right)_*| |\left( W_{l+1}F^{(\alpha)} \right)'|^2
                     dv dv_* d\sigma \right\}^{1/2} \nonumber\\
  & +        \left\{ \iiint_{|v-v_*| \leq 1} |v-v_*|^{\gamma+1} b(cos\theta) \theta^2
                     |\left( F^{(\alpha')} \right)_*| |\left( W_{l-1} F^{(\alpha'')} \right)|^2
                     dv dv_* d\sigma \right\}^{1/2} \nonumber\\
  & ~ \times \left\{ \iiint_{|v-v_*| \leq 1} |v-v_*|^{\gamma+1} b(cos\theta) \theta^2
                     |\left( F^{(\alpha')} \right)_*| |\left( W_{l+1}F^{(\alpha)} \right)'|^2
                     dv dv_* d\sigma \right\}^{1/2} \nonumber\\
\triangleq & G_9^{1/2} \times G_{10}^{1/2} + G_{11}^{1/2} \times G_{12}^{1/2}.\nonumber
\end{align}

It is obvious that
\begin{align*}
  & G_9 \lesssim \| F^{(\alpha')} \|_{L^1} \|W_{l-1}F^{(\alpha'')}\|^2_{L^2},\\
  & G_{10} \lesssim \| F^{(\alpha')} \|_{L^1} \|W_{l+1}F^{(\alpha)}\|^2_{L^2}.
\end{align*}

By Lemma \ref{Lemma 2.6} with $\gamma+1 \in (-3/2,0)$, we have
\begin{align*}
  & G_{11} \lesssim \| F^{(\alpha')} \|_{L^2} \|W_{l-1}F^{(\alpha'')}\|^2_{L^2},\\
  & G_{12} \lesssim \| F^{(\alpha')} \|_{L^2} \|W_{l+1}F^{(\alpha)}\|^2_{L^2}.
\end{align*}

Thus we obtain for $\gamma+1 < 0$,
\begin{align}
  \Psi_{2,1}^{(\alpha',\alpha'')}(t)
\lesssim & \left( \| F^{(\alpha')} \|_{L^1} + \| F^{(\alpha')} \|_{L^2} \right)
           \|W_{l-1} F^{(\alpha'')} \|_{L^2}
           \|W_{l+1} F^{(\alpha)} \|_{L^2} \\
\lesssim & \|W_l F^{(\alpha')} \|_{L^2} \|W_{l-1} F^{(\alpha'')} \|_{L^2}
           \|W_{l+1} F^{(\alpha)} \|_{L^2}. \nonumber
\end{align}

Together with the above estimates for two cases, by setting
$\beta = (1-\gamma^+) {\sf 1}_{\gamma+1 \geq 0} + {\sf 1}_{\gamma+1 < 0}$, we can conclude that, for $\gamma \in (-5/2,1)$,
\begin{align}
  \Psi_{2,1}^{(\alpha',\alpha'')}(t)
\lesssim  \|W_l F^{(\alpha')} \|_{L^2} \|W_{l}F^{(\alpha'')} \|^{\beta}_{L^2}
          \|W_{l+1}F^{(\alpha'')} \|^{1-\beta}_{L^2} \|W_{l+1}F^{(\alpha)} \|_{L^2}. \nonumber
\end{align}

Now we get the estimate of $\Psi_2^{(\alpha',\alpha'')}$ for $\gamma \in (-5/2,1)$,
\begin{align}\label{Psi_2-original}
  \left|\Psi_2^{(\alpha',\alpha'')}(t)\right|
\leq & |\Psi_{2,1}^{(\alpha',\alpha'')}(t)| + |\Psi_{2,2}^{(\alpha',\alpha'')}(t)| \\
\lesssim & \|W_l F^{(\alpha')} \|_{L^2} \|W_{l}F^{(\alpha'')} \|^{\beta}_{L^2}
           \|W_{l+1}F^{(\alpha'')} \|^{1-\beta}_{L^2} \|W_{l+1}F^{(\alpha)} \|_{L^2}, \nonumber
\end{align}
where $\beta=(1-\gamma) \cdot {\sf 1}_{\gamma \in [0,1)} + {\sf 1}_{\gamma \in (-5/2,0)} \in (0,1]$.

We mention again the above inequalities hold for appropriately large $l \geq 4$.

Then we obtain,
\begin{multline}\label{Psi_2}
   \frac{\rho^{2|\alpha|}\big| \Psi_2^{(\alpha',\alpha'')}(t)\big|}{\{(\alpha-r)!\}^{2\nu}}
\leq  C \frac{\{(\alpha'-r)!\}^\nu\{(\alpha''-r)!\}^\nu}{\{(\alpha-r)!\}^\nu}
         \|f(t)\|_{\delta-\kappa t,l,\rho,\alpha',r} \\
         \|f(t)\|^\beta_{\delta-\kappa t,l,\rho,\alpha'',r}
         \|f(t)\|^{1-\beta}_{\delta-\kappa t,l+1,\rho,\alpha'',r}
         \|f(t)\|_{\delta-\kappa t,l+1,\rho,\alpha,r}.
\end{multline}

\begin{rema}
  If we consider only the mild singularity case $s<1/2$, it suffices to obtain the above result by using Taylor expansion of order 1.
\end{rema}

\subsection{The coercivity and commutator estimate}

Considering the estimate on $\Psi_1^{(0,\alpha)}$, we write firstly that
\begin{align*}
\Psi_1^{(0,\alpha)}=&\left\langle Q(\mu F,W_l F^{(\alpha)}),W_l F^{(\alpha)} \right\rangle
                +\left\langle W_l Q(\mu F,F^{(\alpha)})-Q(\mu F,W_l F^{(\alpha)}),W_l F^{(\alpha)} \right\rangle\\
                   =&\Psi_{1,1}^{(0,\alpha)}+\Psi_{1,2}^{(0,\alpha)}.
\end{align*}

\vspace*{1em}
{\noindent\it Step I: The coercivity estimate. }
In order to estimate the ``coercivity term" $\Psi_{1,1}^{(0,\alpha)}$ for soft potential ($\gamma < 0$), using Lemma \ref{coercivity lemma} with $g=\mu F$, $f=W_l F^{(\alpha)}$, we have
\begin{align}
\Psi_{1,1}^{(0,\alpha)}+c_0\| W_{l+\gamma/2} F^{(\alpha)} \|^2_{H^s}
\lesssim &  \left( \| \mu W_{|\gamma+2|} F \|_{L^1}+ \| \mu W_{|\gamma+2|} F \|_{L^2} \right)
             \|W_{l+\gamma/2} F^{(\alpha)}\|^2_{L^2} \\
     & + \left( \| \mu W_{|\gamma|} F \|_{L^1} + \| \mu W_{|\gamma|} F \|_{L^2} \right)
             \| W_{l+\gamma/2} F^{(\alpha)} \|^2_{H^\rho}, \nonumber
\end{align}
where $c_0$ is a constant depending only on the bounds of $\|f\|_{L^1_1}$, $ \|f\|_{L\log L}$.

Due to the definition of the smooth Maxwellian decay solution (Def. \ref{Def1}), we have
\begin{align}
  & \| \mu W_{|\gamma+2|} F \|_{L^1} \leq C\|W_l F \|_{L^2}
=C\|f\|_{\delta-\kappa t,l,\rho,0,r} \leq C, \\
  & \| \mu W_{|\gamma+2|} F \|_{L^2} \leq C\|W_l F \|_{L^2}
=C\|f\|_{\delta-\kappa t,l,\rho,0,r} \leq C,\\
  & \| \mu W_{|\gamma|} F \|_{L^1} \leq C\|W_l F \|_{L^2}
=C\|f\|_{\delta-\kappa t,l,\rho,0,r} \leq C,\\
  & \| \mu W_{|\gamma|} F \|_{L^2} \leq C\|W_l F \|_{L^2}
=C\|f\|_{\delta-\kappa t,l,\rho,0,r} \leq C.
\end{align}
Then we arrive at
\begin{align}
\Psi_{1,1}^{(0,\alpha)}+ c_0 \| W_{l+\gamma/2} F^{(\alpha)} \|^2_{H^s}
\leq C \| W_{l+\gamma/2} F^{(\alpha)} \|^2_{L^2} + C \| W_{l+\gamma/2} F^{(\alpha)} \|^2_{H^\rho}.
\end{align}

Thanks to the following interpolation inequality, for $0<\rho <s$ and any $\varepsilon >0$,
$$
\|f\|_{H^\rho}\leq \|f\|^{\rho/s}_{H^s} \|f\|^{(s-\rho)/s}_{L^2}
              \leq \varepsilon \|f\|_{H^s} + \varepsilon^{-\frac{\rho}{s-\rho}} \|f\|_{L^2},
$$
then choosing $ \varepsilon=c_0/(2C)$, we have for $\gamma < 0$,
\begin{align}\label{coercivity}
\Psi_{1,1}^{(0,\alpha)}+ \frac{c_0}{2} \| W_{l+\gamma/2} F^{(\alpha)} \|^2_{H^s}
\leq C \| W_{l+\gamma/2} F^{(\alpha)} \|^2_{L^2}
\leq C \| W_l F^{(\alpha)} \|^2_{L^2}.
\end{align}

\vspace*{1em}
{\noindent\it Step II: The commutator estimate with $s<1/2$. }
To bound the ``commutator term" $\Psi_{1,2}^{(0,\alpha)}$ in the case of soft potential and $s<1/2$, we need the following proposition (see Proposition 2.8 in \cite{FiveGroup-III}):
\begin{prop}
  Let $0<s<1$, $\gamma > max \{-3, -2s-3/2\}$. For any $l \in \mathbb{R}$,
  \begin{align}
    \left| \langle W_l Q(f,g)-Q(f,W_l g), h \rangle \right|
    \lesssim \|f\|_{L^2_{l+ 3/2+ (2s-1)^{^+} +\gamma^{^+} + \varepsilon}}
             \|g\|_{H^{(2s-1+\varepsilon)^{^+}}_{l+ (2s-1)^{^+} +\gamma^{^+}}} \|h\|_{L^2},
  \end{align}
  where we define $\rho^+=max\{\rho,0\}$.
\end{prop}

Using the above proposition with $f=\mu F$, $g=F^{(\alpha)}$, $h=W_l F^{(\alpha)}$, and taking $\varepsilon$ small enough such that $(2s-1+\varepsilon)^+=0$ in view of $s \in (0,1/2)$, then we have
\begin{align}\label{commutator}
|\Psi_{1,2}^{(0,\alpha)}|
\lesssim  \| \mu W_{l+3/2+\varepsilon} F \|_{L^2} \|W_l F^{(\alpha)} \|^2_{L^2}
\lesssim  \|W_l F^{(\alpha)} \|^2_{L^2} .
\end{align}

Together with (\ref{coercivity}) and (\ref{commutator}) we obtain
\begin{align}
\Psi_1^{(0,\alpha)}+ \frac{c_0}{2} \| W_{l+\gamma/2} F^{(\alpha)} \|^2_{H^s}
\lesssim  \| W_l F^{(\alpha)} \|^2_{L^2} .
\end{align}

Thus we get for $\gamma < 0$,
\begin{align}\label{Psi_1^0}
\frac{\rho^{2|\alpha|}\left| \Psi_1^{(0,\alpha)}(t) \right|}{\{(\alpha-r)!\}^{2\nu}}
+\frac{c_0}{2} \frac{\rho^{2|\alpha|} \| W_{l+\gamma/2} F^{(\alpha)} \|^2_{H^s}}{\{(\alpha-r)!\}^{2\nu}}
\leq C \|f(t)\|^2_{\delta-\kappa t,l,\rho,\alpha,r} .
\end{align}

\vspace*{1em}
{\noindent\it Step III: The commutator estimate with $s=1/2$. }
To bound the ``commutator term" $\Psi_{1,2}^{(0,\alpha)}$ in the critical singularity case $s=1/2$, we restrict our discussion to the condition $-2 \leq \gamma \leq 0$.

We need the following proposition (see Proposition 2.8 in \cite{FiveGroup-III}):
\begin{prop}
  Let $0<s<1$, $\gamma > max \{-3, -2s-3/2\}$. For any $l \in \mathbb{R}$,
  \begin{align}
    \left| \langle W_l Q(f,g)-Q(f,W_l g), h \rangle \right|
    \lesssim \|f\|_{L^2_{l+ 3/2+ (2s-1)^{^+} +\gamma^{^+} + \varepsilon}}
             \|g\|_{H^{(2s-1+\varepsilon)^{^+}}_{l+ (2s-1)^{^+} +\gamma^{^+}}} \|h\|_{L^2},
  \end{align}
  where we define $\rho^+=max\{\rho,0\}$.
\end{prop}

Taking advantage of this proposition with $f=\mu F$, $g=F^{(\alpha)}$, $h=W_l F^{(\alpha)}$, and choosing $\varepsilon=\varepsilon_0$ sufficiently small such that $2\varepsilon_0 < s$, then we have
\begin{align}
|\Psi_{1,2}^{(0,\alpha)}|
\lesssim & \| \mu W_{l+3/2+\varepsilon_{_0}} F \|_{L^2} \|W_l F^{(\alpha)} \|_{H^{\varepsilon_{_0}}_0}
           \|W_l F^{(\alpha)} \|_{L^2}.
\end{align}

Thanks to the interpolation lemma with weight (see Lemma 3.7 in \cite{FiveGroup-Smoothing}):
\begin{lemm}
  For any $k \in \mathbb{R}$, $p \in \mathbb{R}_+$, $\delta > 0$,
  \begin{align*}
    \|f\|^2_{H^k_p(\mathbb{R}^3)} \leq C_\delta \|f\|_{H^{k-\delta}_{2p}(\mathbb{R}^3)}
                                              \|f\|_{H^{k+\delta}_0(\mathbb{R}^3)}.
  \end{align*}
\end{lemm}

Setting $f= W_{l+\gamma/2} F^{(\alpha)}$, $k=\delta=\varepsilon_0$, and $p=-\gamma/2>0$, and using a standard interpolation inequality, we infer that
\begin{align}
  \|W_l F^{(\alpha)} \|^2_{H^{\varepsilon_{_0}}}
\lesssim & \|W_{l+\gamma/2} F^{(\alpha)} \|_{H^{2\varepsilon_{_0}}}
           \|W_{l+\gamma/2} F^{(\alpha)} \|_{L^2_{-\gamma}} \\
\lesssim & \|W_{l+\gamma/2} F^{(\alpha)} \|^{2\varepsilon_0/s}_{H^s}
           \|W_{l+\gamma/2} F^{(\alpha)} \|^{(s-2\varepsilon_0)/s}_{L^2}
           \|W_{l-\gamma/2} F^{(\alpha)} \|_{L^2}. \nonumber
\end{align}

Noticing our assumption $-2 \leq \gamma \leq 0$, this implies
\begin{align}\label{commutator-1/2}
  |\Psi_{1,2}^{(0,\alpha)}|
\lesssim & \|W_l F^{(\alpha)} \|_{L^2} \|W_{l+\gamma/2} F^{(\alpha)} \|^{(s-2\varepsilon_0)/2s}_{L^2}
           \|W_{l-\gamma/2} F^{(\alpha)} \|^{1/2}_{L^2}
           \|W_{l+\gamma/2} F^{(\alpha)} \|^{\varepsilon_0/s}_{H^s} \\
\lesssim & \|W_l F^{(\alpha)} \|^{1+(s-2\varepsilon_0)/2s}_{L^2}
           \|W_{l+1} F^{(\alpha)} \|^{1/2}_{L^2}
           \|W_{l+\gamma/2} F^{(\alpha)} \|^{\varepsilon_0/s}_{H^s}. \nonumber
\end{align}

Together with (\ref{coercivity}) and (\ref{commutator-1/2}), it follows that, for $s=1/2$ and $-2 \leq \gamma \leq 0$,
\begin{align}
  & \Psi_1^{(0,\alpha)}+ \frac{c_0}{2} \| W_{l+\gamma/2} F^{(\alpha)} \|^2_{H^s} \\
\lesssim & \| W_l F^{(\alpha)} \|^2_{L^2}
+ \|W_l F^{(\alpha)} \|^{1+(s-2\varepsilon_0)/2s}_{L^2} \|W_{l+1} F^{(\alpha)} \|^{1/2}_{L^2}
           \|W_{l+\gamma/2} F^{(\alpha)} \|^{\varepsilon_0/s}_{H^s}. \nonumber
\end{align}

Then we get,
\begin{align}\label{Psi_1^0-1/2}
  & \frac{\rho^{2|\alpha|}\left| \Psi_1^{(0,\alpha)}(t) \right|}{\{(\alpha-r)!\}^{2\nu}}
+\frac{c_0}{2}\frac{\rho^{2|\alpha|} \|W_{l+\gamma/2} F^{(\alpha)}\|^2_{H^s}}{\{(\alpha-r)!\}^{2\nu}} \\
\leq & C \|f(t)\|^2_{\delta-\kappa t,l,\rho,\alpha,r}
+ \|f(t)\|^{1+(s-2\varepsilon_0)/2s}_{\delta-\kappa t,l,\rho,\alpha,r}
  \|f(t)\|^{1/2}_{\delta-\kappa t,l+1,\rho,\alpha,r}
  \left( \frac{\rho^{|\alpha|} \|W_{l+\gamma/2} F^{(\alpha)}\|^2_{H^s}}{\{(\alpha-r)!\}^{\nu}} \right)^{\varepsilon_0/s}. \nonumber
\end{align}

\subsection{Upper bound for collision operator}
\vspace*{1em}
{\noindent\it Step I: The upper bound estimate with $s<1/2$. }
In order to estimate $\Psi_{1,1}^{(\alpha',\alpha'')}$ with $s<1/2$, we need the following upper bound estimate for collision operator (see Proposition 2.9 in \cite{FiveGroup-III}):
\begin{prop}\label{upper bound}
  Let $0<s<1$ and $\gamma > max \{-3, -2s-3/2\}$. Then we have, for any $p \in \mathbb{R}$ and $m\in [s-1,s]$,
  \begin{align*}
    \big| \left \langle Q(f,g),h \right \rangle \big|
      \lesssim \left( \|f\|_{L^1_{p^+ + (\gamma+2s)^+}} + \|f\|_{L^2} \right)
               \|g\|_{H^{max{\{s+m,(2s-1+\varepsilon)^+\}}}_{(p+\gamma+2s)^+}} \|h\|_{H^{s-m}_{-p}} .
  \end{align*}
\end{prop}
Using the above proposition with $f=\mu F^{(\alpha')}$, $g=F^{(\alpha'')}$, $h=W_{2l} F^{(\alpha)}$, $p=l-\gamma-2s$, $m=s$, and noticing $\gamma+2s < \gamma +1$, we obtain
\begin{align}
  \left| \Psi_1^{(\alpha',\alpha'')} \right|
 \lesssim &
    \left(\| \mu F^{(\alpha')} \|_{L^1_{l-\gamma -2s+(\gamma+2s)^+}} + \| \mu F^{(\alpha')} \|_{L^2}
    \right) \| W_l F^{(\alpha'')} \|_{H^{2s}}
            \| W_{l+\gamma +2s} F^{(\alpha)} \|_{L^2} \\
 \lesssim & \| W_l F^{(\alpha')} \|_{L^2} \| W_l F^{(\alpha'')} \|_{H^{2s}}
            \| W_{l+1+\gamma} F^{(\alpha)} \|_{L^2},\nonumber
\end{align}
Since $2s<1$ and
$$
\partial_v (W_l \mu^{-1}f^{(\alpha'')})=\partial_v (W_l \mu^{-1}) f^{(\alpha'')}
+ W_l\mu^{-1} f^{(\alpha''+1)},
$$
then we have
$$
\| W_l F^{(\alpha'')} \|_{H^{2s}} \leq \|W_{l+1} F^{(\alpha'')} \|_{L^2} + \|W_l F^{(\alpha''+1)}\|_{L^2}.
$$

Further, by setting $\eta=1+\gamma$ and $\beta=1-\eta^+ \in (0,1]$ as before, we have
$$
\| W_{l+1+\gamma} F^{(\alpha)} \|_{L^2}
\leq \| W_l F^{(\alpha)} \|^\beta_{L^2} \| W_{l+1} F^{(\alpha)} \|^{(1-\beta)}_{L^2}.
$$

So we obtain
\begin{align}
  \left| \Psi_1^{(\alpha',\alpha'')} \right|
 \lesssim & \| W_l F^{(\alpha')} \|_{L^2}  \| W_l F^{(\alpha)} \|^\beta_{L^2}
            \|W_{l+1} F^{(\alpha'')} \|_{L^2} \| W_{l+1} F^{(\alpha)} \|^{(1-\beta)}_{L^2} \\
          & + \| W_l F^{(\alpha')}\|_{L^2} \|W_l F^{(\alpha''+1)}\|_{L^2} \| W_{l+1} F^{(\alpha)}\|_{L^2}
          \nonumber \\
 \triangleq &  \Psi_{1,1}^{(\alpha',\alpha'')}(t) +  \Psi_{1,2}^{(\alpha',\alpha'')}(t).\nonumber
\end{align}

Then we have the estimate for $\Psi_1^{(\alpha',\alpha'')}(t) $ with $\alpha' \neq 0$ as follows:
\begin{multline}\label{Psi_11}
\frac{\rho^{2|\alpha|}\left| \Psi_{1,1}^{(\alpha',\alpha'')}(t) \right|}{\{(\alpha-r)!\}^{2\nu}}
\leq  C \frac{\{(\alpha'-r)!\}^\nu\{(\alpha''-r)!\}^\nu}{\{(\alpha-r)!\}^\nu} \\
 \|f(t)\|_{\delta-\kappa t,l,\rho,\alpha',r}
 \|f(t)\|^{(1-\eta)}_{\delta-\kappa t,l,\rho,\alpha,r}
        \|f(t)\|_{\delta-\kappa t,l+1,\rho,\alpha'',r}
        \|f(t)\|^\eta_{\delta-\kappa t,l+1,\rho,\alpha,r},
\end{multline}
and
\begin{multline}\label{Psi_12}
\frac{\rho^{2|\alpha|} \left| \Psi_{1,2}^{(\alpha',\alpha'')}(t) \right|}{\{(\alpha-r)!\}^{2\nu}}
\leq  C \frac{\{(\alpha'-r)!\}^\nu\{(\alpha''+1-r)!\}^\nu}{\{(\alpha-r)!\}^\nu}\\
 \|f(t)\|_{\delta-\kappa t,l,\rho,\alpha',r} \|f(t)\|_{\delta-\kappa t,l,\rho,\alpha''+1,r}
        \|f(t)\|_{\delta-\kappa t,l+1,\rho,\alpha,r}.
\end{multline}

\vspace*{1em}
{\noindent\it Step II: The upper bound estimate with $s=1/2$. }
Now we limit our discussion to the critical singular case $s=1/2$ and $\gamma+4s \leq 0$.

Using Proposition \ref{upper bound} with $f=\mu F^{(\alpha')}$, $g=F^{(\alpha'')}$, $h=W_{2l} F^{(\alpha)}$, $p=l-\gamma/2$, $m=0$, we obtain
\begin{align}
  \left| \Psi_1^{(\alpha',\alpha'')} \right|
 \lesssim & \left(\| \mu F^{(\alpha')} \|_{L^1_{l-\gamma/2}} + \| \mu F^{(\alpha')} \|_{L^2} \right)
               \| W_l F^{(\alpha'')} \|_{H^s_{\frac{\gamma}{2}+2s}}
               \| W_{l+\gamma/2} F^{(\alpha)} \|_{H^s}  \nonumber\\
 \lesssim & \| W_l F^{(\alpha')} \|_{L^2} \| W_l F^{(\alpha'')} \|_{H^{1/2}}
            \| W_{l+\gamma/2} F^{(\alpha)} \|_{H^s} \nonumber\\
 \lesssim & \| W_l F^{(\alpha')} \|_{L^2} \| W_l F^{(\alpha'')} \|^{1/2}_{H^1}
            \| W_l F^{(\alpha'')} \|^{1/2}_{L^2} \| W_{l+\gamma/2} F^{(\alpha)} \|_{H^s} \nonumber\\
 \lesssim & \| W_l F^{(\alpha')} \|_{L^2} \| W_l F^{(\alpha'')} \|^{1/2}_{L^2}
            \| W_l F^{(\alpha''+1)} \|^{1/2}_{L^2} \| W_{l+\gamma/2} F^{(\alpha)} \|_{H^s} \nonumber\\
        & + \| W_l F^{(\alpha')} \|_{L^2} \| W_l F^{(\alpha'')} \|^{1/2}_{L^2}
            \| W_{l+1} F^{(\alpha'')} \|^{1/2}_{L^2} \| W_{l+\gamma/2} F^{(\alpha)} \|_{H^s} \nonumber\\
 \triangleq &  \Psi_{1,1}^{(\alpha',\alpha'')}(t) +  \Psi_{1,2}^{(\alpha',\alpha'')}(t),\nonumber
\end{align}
where we have used a standard interpolation inequality and
$$
\partial_v (W_l \mu^{-1}f^{(\alpha'')})=\partial_v (W_l \mu^{-1}) f^{(\alpha'')}
+ W_l\mu^{-1} f^{(\alpha''+1)}.
$$

Then we get the estimate for $\Psi_1^{(\alpha',\alpha'')}(t) $ with $\alpha' \neq 0$ as follows:
\begin{multline}\label{Psi_11-1/2}
 \frac{\rho^{2|\alpha|}\left| \Psi_{1,1}^{(\alpha',\alpha'')}(t) \right|}{\{(\alpha-r)!\}^{2\nu}}
\leq  C \frac{\{(\alpha'-r)!\}^\nu\{(\alpha''-r)!\}^\nu}{\{(\alpha-r)!\}^\nu}
\|f(t)\|_{\delta-\kappa t,l,\rho,\alpha',r} \\
\|f(t)\|^{1/2}_{\delta-\kappa t,l,\rho,\alpha'',r}
\|f(t)\|^{1/2}_{\delta-\kappa t,l,\rho,\alpha''+1,r}
        \big( \frac{\rho^{|\alpha|} \|W_{l+\gamma/2} F^{(\alpha)}\|_{H^s}}{\{(\alpha-r)!\}^{\nu}} \big),
\end{multline}
and
\begin{multline}\label{Psi_12-1/2}
 \frac{\rho^{2|\alpha|} \left| \Psi_{1,2}^{(\alpha',\alpha'')}(t) \right|}{\{(\alpha-r)!\}^{2\nu}}
\leq  C \frac{\{(\alpha'-r)!\}^\nu\{(\alpha''+1-r)!\}^\nu}{\{(\alpha-r)!\}^\nu}
    \|f(t)\|_{\delta-\kappa t,l,\rho,\alpha',r} \\
    \|f(t)\|^{1/2}_{\delta-\kappa t,l,\rho,\alpha'',r}
    \|f(t)\|^{1/2}_{\delta-\kappa t,l+1,\rho,\alpha'',r}
        \big( \frac{\rho^{|\alpha|} \|W_{l+\gamma/2} F^{(\alpha)}\|_{H^s}}{\{(\alpha-r)!\}^{\nu}} \big).\
\end{multline}

Now we have completed the estimates for $\Psi_1^{(0,\alpha)}(t)$, $\Psi_1^{(\alpha',\alpha'')}(t)$,
and $\Psi_2^{(\alpha',\alpha'')}(t)$. We remark that the assumptions $\gamma+4s \leq 0$ and $\gamma \in [-2,0]$ give the final assumption $\gamma=-2$ when $s=1/2$. By the way, this implies that $\beta=1$.

\subsection{Completion of the proof of the main lemma}
First of all, we refer to Lemma 3.6 of \cite{Zhang-Yin} (see also, Proposition 3.1 of \cite{Mori-Ukai}):
\begin{prop}\label{Coefficient}
If $\nu \geq 1$ and $2\leq r \in \mathbb{N}$ then there exists a constant $B>0$ depending only on $n$ and $r$ such that for any $\alpha \in \mathbb{Z}^n$
\begin{align}
\sum_{\alpha=\alpha'+\alpha''} \frac{\alpha!}{\alpha'!\alpha''!}
      \frac{\{(\alpha'-r)!\}^\nu\{(\alpha''-r)!\}^\nu}{\{(\alpha-r)!\}^\nu} \leq B,
\end{align}
Furthermore, if $\nu>1$ and $r>1+\nu/(\nu-1)$ then there exists a constant $B'>0$ depending only on $n$, $\nu$ and $r$ such that for any $0 \neq \alpha \in \mathbb{Z}^n$
\begin{align}
\sum_{\alpha=\alpha'+\alpha'',\,\alpha' \neq 0} \frac{\alpha!}{\alpha'!\alpha''!}
      \frac{\{(\alpha'-r)!\}^\nu\{(\alpha''+1-r)!\}^\nu}{\{(\alpha-r)!\}^\nu} \leq B'.
\end{align}
\end{prop}

Now we resume the proof of the main lemma (Lemma \ref{Main Lemma}). Firstly we consider the integral including the term $\mathcal{K}(t)$ on the right-hand side of (\ref{EQ2}), we get from (\ref{Psi_2}) that
\begin{align}\label{K-estimate}
&\int_0^t \frac{\rho^{2|\alpha|} \mathcal{K}(\tau)}{\{(\alpha-r)!\}^{2\nu}}d\tau \\
\leq & C \sum \frac{\alpha!}{\alpha'! \alpha''!}
              \frac{\{(\alpha'-r)!\}^\nu\{(\alpha''-r)!\}^\nu}{\{(\alpha-r)!\}^\nu}
         \int_0^t \|f(\tau)\|_{\delta-\kappa \tau,l,\rho,\alpha',r}
                  \|f(\tau)\|^\beta_{\delta-\kappa \tau,l,\rho,\alpha'',r} \nonumber\\
& \hspace*{15em} \times \|f(\tau)\|^{1-\beta}_{\delta-\kappa \tau,l+1,\rho,\alpha'',r}
                        \|f(\tau)\|_{\delta-\kappa \tau,l+1,\rho,\alpha,r}d\tau \nonumber\\
\leq & CB \int_0^t \|f(\tau)\|_{\delta-\kappa \tau,l,\rho,\alpha',r}
                   \|f(\tau)\|^\beta_{\delta-\kappa \tau,l,\rho,\alpha'',r}
                   \|f(\tau)\|^{1-\beta}_{\delta-\kappa \tau,l+1,\rho,\alpha'',r}
                   \|f(\tau)\|_{\delta-\kappa \tau,l+1,\rho,\alpha,r}d\tau \nonumber\\
\leq & CB \Big\{ \frac{1}{4\varepsilon}
                  \int_0^t \|f(\tau)\|^{2/\beta}_{\delta-\kappa \tau,l,\rho,\alpha',r}
                           \|f(\tau)\|^2_{\delta-\kappa \tau,l,\rho,\alpha'',r} d\tau
 + \varepsilon \int_0^t
                  \left(\|f(\tau)\|^2_{\delta-\kappa \tau,l+1,\rho,\alpha'',r}
                        +\|f(\tau)\|^2_{\delta-\kappa \tau,l+1,\rho,\alpha,r} \right) d\tau
          \Big\}. \nonumber
\end{align}

Secondly we consider the integral including the term $\mathcal{J}(t)$ on the right-hand side of (\ref{EQ2}), we get from (\ref{Psi_11}) and (\ref{Psi_12}) that, for $s<1/2$,
\begin{align}\label{J-estimate}
& \int_0^t \frac{\rho^{2|\alpha|} \mathcal{J}(\tau)}{\{(\alpha-r)!\}^{2\nu}}d\tau \\
\leq & \int_0^t \sum_{\alpha=\alpha'+\alpha''} \frac{\alpha!}{\alpha'!\alpha''!}
        \frac{\rho^{2|\alpha|}
        \left| \Psi_{1,1}^{(\alpha',\alpha'')}(\tau) \right|}{\{(\alpha-r)!\}^{2\nu}}d\tau
+ \int_0^t \sum_{\alpha=\alpha'+\alpha'',\,\alpha' \neq 0}
        \frac{\alpha!}{\alpha'!\alpha''!}\frac{\rho^{2|\alpha|}
        \left| \Psi_{1,2}^{(\alpha',\alpha'')}(\tau) \right|}{\{(\alpha-r)!\}^{2\nu}}d\tau \nonumber\\
\leq & CB \int_0^t \|f(\tau)\|_{\delta-\kappa t,l,\rho,\alpha',r}
                   \|f(\tau)\|^\beta_{\delta-\kappa t,l,\rho,\alpha,r}
                   \|f(\tau)\|_{\delta-\kappa t,l+1,\rho,\alpha'',r}
                   \|f(\tau)\|^{(1-\beta)}_{\delta-\kappa t,l+1,\rho,\alpha,r} d\tau \nonumber\\
  & + CB' \int_0^t \|f(\tau)\|_{\delta-\kappa t,l,\rho,\alpha',r}
                   \|f(\tau)\|_{\delta-\kappa t,l,\rho,\alpha''+1,r}
                   \|f(\tau)\|_{\delta-\kappa t,l+1,\rho,\alpha,r} d\tau\nonumber\\
\leq & C \Big\{ \frac{1}{4\varepsilon}
                  \int_0^t \|f(\tau)\|^{2/\beta}_{\delta-\kappa \tau,l,\rho,\alpha',r}
                           \|f(\tau)\|^2_{\delta-\kappa \tau,l,\rho,\alpha'',r} d\tau
 + \varepsilon \int_0^t
                  \left(\|f(\tau)\|^2_{\delta-\kappa \tau,l+1,\rho,\alpha'',r}
                        +\|f(\tau)\|^2_{\delta-\kappa \tau,l+1,\rho,\alpha,r} \right) d\tau
          \Big\} \nonumber\\
& + C \Big\{ \frac{1}{4\varepsilon}
                  \int_0^t \|f(\tau)\|^2_{\delta-\kappa \tau,l,\rho,\alpha',r}
                           \|f(\tau)\|^2_{\delta-\kappa \tau,l,\rho,\alpha''+1,r} d\tau
 + \varepsilon \int_0^t
                  \|f(\tau)\|^2_{\delta-\kappa \tau,l+1,\rho,\alpha,r}  d\tau
          \Big\}.\nonumber 
\end{align}

And from (\ref{Psi_11-1/2}) and (\ref{Psi_12-1/2}), we infer that, for $s=1/2$,
\begin{align}
& \int_0^t \frac{\rho^{2|\alpha|} \mathcal{J}(\tau)}{\{(\alpha-r)!\}^{2\nu}}d\tau \\
\leq & \int_0^t \sum_{\alpha=\alpha'+\alpha'',~\alpha' \neq 0} \frac{\alpha!}{\alpha'!\alpha''!}
        \frac{\rho^{2|\alpha|}
        \left| \Psi_{1,1}^{(\alpha',\alpha'')}(\tau) \right|}{\{(\alpha-r)!\}^{2\nu}}d\tau
+ \int_0^t \sum_{\alpha=\alpha'+\alpha'',~\alpha' \neq 0}
        \frac{\alpha!}{\alpha'!\alpha''!}\frac{\rho^{2|\alpha|}
        \left| \Psi_{1,2}^{(\alpha',\alpha'')}(\tau) \right|}{\{(\alpha-r)!\}^{2\nu}}d\tau \nonumber\\
\leq & CB' \int_0^t \|f(\tau)\|_{\delta-\kappa \tau,l,\rho,\alpha',r}
                   \|f(\tau)\|^{1/2}_{\delta-\kappa \tau,l,\rho,\alpha'',r}
                   \|f(\tau)\|^{1/2}_{\delta-\kappa \tau,l,\rho,\alpha''+1,r}
                   \left( \frac{\rho^{|\alpha|} \|W_{l+\gamma/2} F^{(\alpha)}\|_{H^s}}{\{(\alpha-r)!\}^{\nu}} \right) d\tau \nonumber\\
  & + CB' \int_0^t \|f(\tau)\|_{\delta-\kappa \tau,l,\rho,\alpha',r}
                   \|f(\tau)\|^{1/2}_{\delta-\kappa \tau,l,\rho,\alpha'',r}
                   \|f(\tau)\|^{1/2}_{\delta-\kappa \tau,l+1,\rho,\alpha'',r}
                   \left( \frac{\rho^{|\alpha|} \|W_{l+\gamma/2} F^{(\alpha)}\|_{H^s}}{\{(\alpha-r)!\}^{\nu}} \right) d\tau \nonumber\\
\triangleq & \mathcal{J}_1(t) + \mathcal{J}_2(t). \nonumber
\end{align}

By the H\"{o}lder inequality we have
\begin{align}\label{J1-estimate-1/2}
  \mathcal{J}_1(t)
\leq CB' \Big\{ & C_\varepsilon
                  \int_0^t \left(\|f(\tau)\|^4_{\delta-\kappa \tau,l,\rho,\alpha',r}
                           \|f(\tau)\|^2_{\delta-\kappa \tau,l,\rho,\alpha'',r}
                        +  \|f(\tau)\|^2_{\delta-\kappa \tau,l,\rho,\alpha''+1,r} \right)d\tau \\
              & +  \varepsilon
              \int_0^t \frac{\rho^{2|\alpha|} \|W_{l+\gamma/2} F^{(\alpha)}\|^2_{H^s}}
                                        {\{(\alpha-r)!\}^{2\nu}} d\tau
        \Big\}, \nonumber
\end{align}
and
\begin{align}\label{J2-estimate-1/2}
  \mathcal{J}_2(t)
\leq CB' \Big\{ & C_\varepsilon
                  \int_0^t \left( \|f(\tau)\|^4_{\delta-\kappa \tau,l,\rho,\alpha',r}
                                  \|f(\tau)\|^2_{\delta-\kappa \tau,l,\rho,\alpha'',r}
                           \right) d\tau \\
      & + \varepsilon \left( \int_0^t \|f(\tau)\|^2_{\delta-\kappa \tau,l+1,\rho,\alpha'',r}d\tau
                          + \int_0^t \frac{\rho^{2|\alpha|} \|W_{l+\gamma/2} F^{(\alpha)}\|^2_{H^s}}
                                          {\{(\alpha-r)!\}^{2\nu}} d\tau
                      \right)
        \Big\}. \nonumber
\end{align}

Since (\ref{Psi_1^0-1/2}) implies
\begin{align}\label{Psi_1^0-estimate-1/2}
  & \int_0^t \frac{\rho^{2|\alpha|}\left| \Psi_1^{(0,\alpha)}(\tau) \right|}{\{(\alpha-r)!\}^{2\nu}} d\tau
+\frac{c_0}{2} \int_0^t \frac{\rho^{2|\alpha|} \|W_{l+\gamma/2} F^{(\alpha)}\|^2_{H^s}}{\{(\alpha-r)!\}^{2\nu}} d\tau\\
\leq & C \int_0^t \|f(\tau)\|^2_{\delta-\kappa \tau,l,\rho,\alpha,r} d\tau
 + C \int_0^t \|f(\tau)\|^{1+(s-2\varepsilon_0)/2s}_{\delta-\kappa \tau,l,\rho,\alpha,r}
  \|f(\tau)\|^{1/2}_{\delta-\kappa \tau,l+1,\rho,\alpha,r}
  \left( \frac{\rho^{|\alpha|} \|W_{l+\gamma/2} F^{(\alpha)}\|^2_{H^s}}{\{(\alpha-r)!\}^{\nu}} \right)^{\varepsilon_0/s} d\tau \nonumber\\
\leq & C_\varepsilon \int_0^t \|f(\tau)\|^2_{\delta-\kappa \tau,l,\rho,\alpha,r} d\tau
 + \varepsilon \left( \int_0^t \|f(\tau)\|^2_{\delta-\kappa \tau,l+1,\rho,\alpha,r}d\tau
                          + \int_0^t \frac{\rho^{2|\alpha|} \|W_{l+\gamma/2} F^{(\alpha)}\|^2_{H^s}}
                                          {\{(\alpha-r)!\}^{2\nu}} d\tau
                      \right). \nonumber
\end{align}

Combining (\ref{EQ2}) with this inequality, (\ref{K-estimate}), (\ref{J1-estimate-1/2}), and (\ref{J2-estimate-1/2}), we obtain, with $s=1/2$ and a small enough $\varepsilon$:
\begin{align}\label{Inequality-1/2}
&\|f(t)\|^2_{\delta-\kappa t,l,\rho,\alpha,r}
+ \frac{c_0}{4} \int_0^t
      \frac{\rho^{2|\alpha|} \| W_{l+\gamma/2} F^{(\alpha)} \|^2_{H^s}}{\{(\alpha-r)!\}^{2\nu}}
+ 2 \kappa \int_0^t \|f(\tau)\|^2_{\delta-\kappa \tau,l+1,\rho,\alpha,r} d\tau \\
\leq &  \|f(0)\|^2_{\delta,l,\rho,\alpha,r}
   + C_\varepsilon
   \int_0^t \Big\{ \|f(\tau)\|^{2/\beta}_{\delta-\kappa t,l,\rho,\alpha',r}
                           \|f(\tau)\|^2_{\delta-\kappa t,l,\rho,\alpha'',r} \nonumber\\
     & + \|f(\tau)\|^4_{\delta-\kappa t,l,\rho,\alpha',r}
                             \|f(\tau)\|^2_{\delta-\kappa t,l,\rho,\alpha'',r}
                           +  \|f(\tau)\|^2_{\delta-\kappa t,l,\rho,\alpha''+1,r}
     + \|f(\tau)\|^2_{\delta-\kappa t,l,\rho,\alpha,r} \Big\} d\tau
   \nonumber\\
 & + \varepsilon
   \int_0^t \Big\{ \|f(\tau)\|^2_{\delta-\kappa t,l+1,\rho,\alpha'',r}
                        +\|f(\tau)\|^2_{\delta-\kappa t,l+1,\rho,\alpha,r} \Big\} d\tau, \nonumber
\end{align}

On the other hand, combining \eqref{EQ2} with \eqref{Psi_1^0}, \eqref{K-estimate} and \eqref{J-estimate} gives that, for $s<1/2$,
\begin{align}\label{Inequality}
&\|f(t)\|^2_{\delta-\kappa t,l,\rho,\alpha,r}
+ \frac{c_0}{2} \int_0^t
      \frac{\rho^{2|\alpha|} \| W_{l+\gamma/2} F^{(\alpha)} \|^2_{H^s}}{\{(\alpha-r)!\}^{2\nu}}
+ 2 \kappa \int_0^t \|f(\tau)\|^2_{\delta-\kappa \tau,l+1,\rho,\alpha,r} d\tau \\
\leq &  \|f(0)\|^2_{\delta,l,\rho,\alpha,r}
   + C_\varepsilon
       \int_0^t \Big\{ \|f(\tau)\|^{2/\beta}_{\delta-\kappa t,l,\rho,\alpha',r}
                           \|f(\tau)\|^2_{\delta-\kappa t,l,\rho,\alpha'',r}
                     + \|f(\tau)\|^2_{\delta-\kappa t,l,\rho,\alpha',r}
                             \|f(\tau)\|^2_{\delta-\kappa t,l,\rho,\alpha''+1,r} \nonumber\\
                  &  + \|f(\tau)\|^2_{\delta-\kappa t,l,\rho,\alpha,r}
                  \Big\}d\tau
 + \varepsilon
   \int_0^t \Big\{ \|f(\tau)\|^2_{\delta-\kappa t,l+1,\rho,\alpha'',r}
                        +\|f(\tau)\|^2_{\delta-\kappa t,l+1,\rho,\alpha,r} \Big\} d\tau, \nonumber
\end{align}

Noticing that $3r \leq |\alpha|=|\alpha'|+|\alpha''| \leq N $, letting
$$
A=A_r(f)=\sup_{t\in [0,T]} \max_{|\beta| \leq 3r} \|f(t)\|_{\delta-\kappa t,l,\rho,\beta,r},
$$
we will give a different estimate on the factor
$$
\Xi \triangleq \|f(\tau)\|^{2/\beta}_{\delta-\kappa t,l,\rho,\alpha',r}
 \|f(\tau)\|^2_{\delta-\kappa t,l,\rho,\alpha'',r},
$$
with respect to $|\alpha'|$ and $|\alpha''|$ while taking supremum on $|\alpha|$:
\begin{itemize}
  \item[$\cdot$] if $|\alpha'|,~|\alpha''|<3r$, we have
    $\underset{3r \leq |\alpha| \leq N}{\sup} \Xi \leq A^{2(1+\beta)/\beta};$
  \item[$\cdot$] if $|\alpha'|<3r,~|\alpha''| \geq 3r$, we have
    $\underset{3r\leq|\alpha| \leq N}{\sup}\Xi\leq
        A^{2/\beta}\|f\|^2_{l,\rho,r,N}(\tau);$
  \item[$\cdot$] if $|\alpha'| \geq 3r,~|\alpha''| < 3r$, we have
    $\underset{3r \leq |\alpha| \leq N}{\sup} \Xi \leq A^2
        \|f\|^{2/\beta}_{l,\rho,r,N}(\tau);$
  \item[$\cdot$] if $|\alpha'|,~|\alpha''| \geq 3r$, we have
    $\underset{3r \leq |\alpha| \leq N}{\sup} \Xi \leq \|f\|^{2(1+\beta)/\beta}_{l,\rho,r,N}(\tau).$
\end{itemize}

Then we obtain
\begin{align}
  \underset{3r \leq |\alpha| \leq N}{\sup} \Xi
  \leq \|f\|^2_{l,\rho,r,N}(\tau) + \|f\|^{2(1+\beta)/\beta}_{l,\rho,r,N}(\tau).
\end{align}

Other terms $\|f(\tau)\|^4_{\delta-\kappa t,l,\rho,\alpha',r} \|f(\tau)\|^2_{\delta-\kappa t,l,\rho,\alpha'',r}$ and $\|f(\tau)\|^2_{\delta-\kappa t,l,\rho,\alpha',r} \|f(\tau)\|^2_{\delta-\kappa t,l,\rho,\alpha''+1,r}$ can be controlled similarly, inserting these results into \eqref{Inequality} and \eqref{Inequality-1/2}, we finally obtain the Gronwall type inequality:
\begin{align}
&\|f(t)\|^2_{\delta-\kappa t,l,\rho,\alpha,r}
+ \frac{c_0}{4} \int_0^t
      \frac{\rho^{2|\alpha|} \| W_{l+\gamma/2} \mu^{-1} f^{(\alpha)}(\tau) \|^2_{H^s}}
           {\{(\alpha-r)!\}^{2\nu}} d\tau
+ 2 \kappa \int_0^t \|f(\tau)\|^2_{\delta-\kappa \tau,l+1,\rho,\alpha,r} d\tau \\
\leq &  \|f(0)\|^2_{\delta,l,\rho,\alpha,r}
       + C_\kappa \int_0^t \left(\|f\|^2_{l,\rho,r,N}(\tau) + \|f\|^{2(i+\beta)/\beta}_{l,\rho,r,N}(\tau)
                          \right) d\tau
       + \frac{\kappa}{10}\sup_{3r \leq |\alpha| \leq N}
                  \int_0^t \|f(\tau)\|^2_{\delta-\kappa \tau,l+1,\rho,\alpha,r} d\tau, \nonumber
\end{align}
where $i=1$ when $s<1/2$ or $i=2$ if $s=1/2$. This leads to the desired estimate (\ref{basic inequa}) including the extra second term of the left-hand side.

This completes the proof of Lemma \ref{Main Lemma}.\qed

\bigskip
\noindent\textbf{Acknowledgements.}  This work was partially supported by NNSFC (No. 11271382 and No. 10971235), RFDP (No. 20120171110014), and the key project of Sun Yat-sen University (No. c1185).

\phantomsection

\end{document}